\documentclass[10pt]{article}
 
\usepackage[latin1]{inputenc}
\usepackage[title]{appendix}
\usepackage{url, amsmath,enumerate,fancyhdr,amssymb, amsthm, 
pstricks, epsf, lscape, ifpdf, graphicx, float}
\usepackage{paralist}
\usepackage[margin=3.0cm]{geometry} 

\usepackage{algorithm}
\usepackage{algpseudocode}
\algnewcommand\algorithmicinput{\textbf{Input:}}
\algnewcommand\INPUT{\item[\algorithmicinput]}
\algnewcommand\algorithmicinitialization{\textbf{Initialization:}}
\algnewcommand\INITIALIZATION{\item[\algorithmicinitialization]}

\renewcommand{\arraystretch}{1.2}

\newcommand{\compilepdf}[1]{#1} 

\mathchardef\mhyphen="2D

\newtheorem{Definition}{Definition}

\newtheorem{Example}{Example}
\newtheorem{Proposition}{Proposition}
\newtheorem{Lemma}{Lemma}
\newtheorem{Theorem}{Theorem}
\newtheorem{Corollary}{Corollary}
\newtheorem{Remark}{Remark}
\newtheorem{Assumption}{Assumption}

\newcommand{\palt}{P_{\rm alt}}
\newcommand{\dalt}{D_{\rm alt}}

\newcommand{\un}{\underline} 
\newcommand{\lin}{\operatorname{lin}}
\newcommand{\primal}{\operatorname{Primal}} 

\newcommand{\front}{\operatorname{front}}
\newcommand{\codim}{\operatorname{codim}}
\newcommand{\trace}{\operatorname{trace}}
\newcommand{\fr}{\operatorname{FR}}

\newcommand{\cl}{\operatorname{cl}}
\newcommand{\ri}{\operatorname{ri}}
\newcommand{\dir}{\operatorname{dir}}

\newcommand{\eps}{\epsilon} 
\newcommand{\bpx}{\begin{pmatrix}}
\newcommand{\epx}{\end{pmatrix}}
\newcommand{\bbx}{\begin{bmatrix}}
\newcommand{\ebx}{\end{bmatrix}}

\newcommand{\aut}{\operatorname{Aut}}

\newcommand{\bdef}{\begin{Definition}} 
\newcommand{\commentout}[1]{}
\newcommand{\co}[1]{}

\newcommand{\nin}{\noindent}
\newcommand{\ti}{\times}
 
\newcommand{\pf}[1]{\vspace{.35cm} \nin {\bf Proof {#1} }}

\newcommand{\norm}[1]{\parallel \! #1 \! \parallel}
\newcommand{\sym}[1]{{\cal S}^{#1}}
\newcommand{\psd}[1]{{\cal S}_+^{#1}}
\newcommand{\rad}[1]{\mathbb{R}^{#1}}
\newcommand{\zad}[1]{\mathbb{Z}^{#1}}
\newcommand{\so}[1]{\mathbb{SO}({#1})}
\newcommand{\eref}[1]{(\ref{#1})}

\newcommand{\R}{ {\cal R} }
\newcommand{\N}{ {\cal N} }

\newcommand{\la}{\langle}
\newcommand{\ra}{\rangle}

\newcommand{\by}{\bar{y}}

\newcommand{\tri}{\triangle}

\newcommand{\beq}{\begin{equation}}
\newcommand{\eeq}{\end{equation}}
\newcommand{\beqa}{\begin{eqnarray}}
\newcommand{\eeqa}{\end{eqnarray}}
\newcommand{\ba}{\begin{array}}
\newcommand{\ena}{\end{array}}
\newcommand{\bac}{\begin{array}{ccccccccccc}}
\newcommand{\eac}{\end{array}}
\newcommand{\bprop}{\begin{Proposition}}
\newcommand{\eprop}{\end{Proposition}}

\newcommand{\beqast}{\begin{eqnarray*}}
\newcommand{\eeqast}{\end{eqnarray*}}
\newcommand{\benum}{\begin{enumerate}}
\newcommand{\eenum}{\end{enumerate}}
\newcommand{\bit}{\begin{itemize}}
\newcommand{\eit}{\end{itemize}}
\newcommand{\bth}{\begin{Theorem}}
\newcommand{\enth}{\end{Theorem}}
\newcommand{\ble}{\begin{Lemma}}
\newcommand{\ele}{\end{Lemma}}
\newcommand{\bex}{\begin{Example}}
\newcommand{\eex}{\end{Example}}
\newcommand{\bcor}{\begin{Corollary}}
\newcommand{\ecor}{\end{Corollary}}
\newcommand{\brem}{\begin{Remark}}
\newcommand{\erem}{\end{Remark}}
\newcommand{\bass}{\begin{Assumption}}
\newcommand{\eass}{\end{Assumption}}

\newcommand{\LRA}{\Leftrightarrow}

\renewcommand{\arraystretch}{1.2}

\newcommand\emp{\emptyset}

\newcommand\clps{\hspace{3.5cm}}

\newcommand{\bsmx}{\begin{small} \begin{pmatrix}}
\newcommand{\esmx}{\end{pmatrix} \end{small}}
\newcommand{\regfr}{\operatorname{REGFR}}
\newcommand{\relregfr}{\operatorname{RELREGFR}}
\newcommand{\revregfr}{\operatorname{REVREGFR}}

\title{\Large Exact duals and 
short certificates of infeasibility and weak infeasibility in conic linear programming}

\author{Minghui Liu \thanks{liu.m.h2010@gmail.com, SAS Inc., Cary} \hspace{2.5cm} G\'{a}bor Pataki  \thanks{gabor@unc.edu, Department of Statistics and Operations Research, 
		University of North Carolina at Chapel Hill }  }

\begin{document}

\co{Change vs april 10: add section on generation}

\newcounter{algorithmcounter}

\maketitle 

\begin{abstract} 
	
		In conic linear programming -- in contrast to linear programming -- 
	the Lagrange dual is not an exact dual: it may not attain its optimal value, or there may be a positive duality gap. 
	The corresponding Farkas' lemma is also not exact (it does not always prove infeasibility).
	\co{
		In conic linear programming -- in contrast to linear programming -- the Lagrange dual is not an exact dual (it may not attain its optimal value, or a positive duality gap may occur), and the corresponding Farkas' lemma may fail to prove infeasibility. We describe exact duals, and certificates of infeasibility and weak infeasibility for conic LPs which are nearly as simple as the Lagrange dual, but do not rely on any constraint qualification. 
			}
		We describe exact duals, and exact certificates of infeasibility and weak infeasibility for conic LPs which are nearly as simple as the Lagrange dual, but do not rely on any constraint qualification. 
	Some of our exact duals generalize the SDP duals of Ramana, Klep and Schweighofer to the 
	context of general conic LPs. 
	Some of our infeasibility certificates generalize 
	the row echelon form of a linear system of equations: they consist of a small, trivially 
	infeasible subsystem obtained by elementary row operations. 
		We prove analogous results for weakly infeasible systems. 
	
	We obtain some fundamental geometric corollaries: an exact characterization of when the 
	linear image of a closed convex cone is closed, and an exact characterization of nice cones.
	
	Our infeasibility certificates provide algorithms 
	to generate {\em all} infeasible conic LPs over several important classes of cones; and 
	{\em all} weakly infeasible SDPs in a natural class. 
	Using these algorithms 
	we generate a public domain library of infeasible and 
	weakly infeasible SDPs. The status of our instances can be verified 
     by inspection in exact arithmetic, 
	but they turn out to be challenging for commercial and research codes.
	
\end{abstract}



{\em OR/MS subject classification:} Primary: convexity; secondary: programming-nonlinear-theory


\section{Introduction and a sample of the main results}

Many problems in engineering, combinatorial optimization, and economics 
can be expressed as a conic linear program (LP) of the form 
\beq \label{p} \tag{P} 
\ba{rl} 
\sup & \la c, x\ra \\
s.t. & Ax \leq_{K} b, \\
\ena
\eeq
where $A: \rad{m} \rightarrow Y$ is a linear map, 
$Y$ is a finite dimensional Euclidean space, $K \subseteq Y$ is a closed convex cone, and 
$s \leq_K t$ stands for $t - s \in K.$ 
We naturally associate a dual program with \eref{p}: 
letting $A^*$ be the adjoint of $A, \,$ and 
$K^*$ the dual cone of $K, \,$ its Lagrange dual is 
\beq \label{d} \tag{D} 
\ba{rrcl} 
\inf & \la b, y\ra \\
s.t. &  A^{*}y &=& c \\
     & y & \geq_{K^*} & 0.
\ena
\eeq
Problems \eref{p} and \eref{d} generalize linear programs, and 
weak duality -- the inequality $\la c, x \ra \leq \la b, y \ra$ between 
a pair of feasible solutions -- trivially holds. 
However, in contrast to linear programming,
the optimal values of \eref{p} and of \eref{d} may differ, and/or may not be attained.

A suitable conic linear system can prove the infeasibility of \eref{p} and of 
\eref{d}: their classical alternative systems are
$$
\ba{rrclcrcll}
&   A^* y  & =  & 0 & \clps  & A x  & \geq_K & 0 & \\
(\palt)        &    \la b, y \ra & =  & -1 & {\rm and}   & \la c, x \ra & = & -1. & (\dalt) \\
&   y        & \geq_{K^*} & 0 & 
\ena
$$
When $(\palt)$ is feasible, \eref{p} is trivially infeasible, and we call it
{\em strongly infeasible.} 
However -- again in contrast to linear programming -- 
$(\palt)$ and \eref{p} may both be infeasible, 
and in this case we call \eref{p} {\em weakly infeasible}. We define strong and weak infeasibility of \eref{d} analogously.

These pathological behaviors -- nonattainment of the optimal values, positive duality gaps, and weak infeasibility -- occur in semidefinite programs (SDPs) and second order conic programs, which are arguably the most useful classes of conic linear programs. Pathological conic LPs are often  difficult, or impossible to solve.

This paper focuses on exact duals and exact certificates of infeasibility and weak infeasibility of \eref{p} and 
of \eref{d}. An {\em exact dual, or strong dual} of \eref{p} is a conic LP with an ''inf'' objective,
which i) satisfies weak duality, ii) has the same optimal value as \eref{p}, and iii) attains this value, when it is finite. 
We define an exact dual of \eref{d} analogously.

An {\em exact certificate of infeasibility} of a conic LP 
is a finite set of vectors, from which a 
suitable polynomial time algorithm (a "verifier") 
can deduce that the conic LP is indeed infeasible.
The term "exact" stresses that 
such a set of vectors must exist for {\em every} infeasible instance. 
We  define 
exact certificates 
of other properties (say of weak infeasibility) of conic LPs analogously. 
The latter definitions are  informal, but they  will be enough 
to explain our ideas and results;   we also 
give a formal definition in Appendix \ref{app-def}. For simplicity, we will often talk about the 
infeasibility of a conic linear {\em system.}

Exact certificates of infeasibility of conic LPs, say, of \eref{p} will appear in one of the 
following forms in this paper:
\begin{itemize}
	\item Either as a conic linear system which is feasible exactly when \eref{p} is infeasible;
	\item Or as a transformation of \eref{p} into an equivalent problem, whose infeasibility is  "easy" to verify.
\end{itemize}
Exact duals and exact certificates of infeasibility of conic LPs are, 
of course,  known in special cases. 
For example, if $K$ is polyhedral,  then \eref{d} is an exact dual of \eref{p}, 
$(\palt)$ is an exact certificate of infeasibility of \eref{p}, and 
$(\dalt)$ is an exact certificate of infeasibility of \eref{d}. 
If $K = \{0 \},$  then \eref{p} is a linear system of equations; then  
both $(\palt),$ and the row echelon form of \eref{p} (which contains an obviously infeasible equation $\la 0, x \ra = 1$) are exact certificates of infeasibility. 
However,  in general,  \eref{d} is not an exact dual,  and $(\palt)$ and $(\dalt)$ are  not exact certificates of
infeasibility. Moreover, general conic LPs are not known to have an equivalent of 
a row echelon form.

 This paper builds on three approaches, which 
 provide exact duals, and exact certificates of infeasibility for conic LPs
 (and which we review below, with other relevant 
 references): the first is
facial reduction algorithms -- see Borwein and Wolkowicz \cite{BorWolk:81,BorWolk:81B}, Waki and Muramatsu \cite{WakiMura:12}, Pataki \cite{Pataki:13}; and the second is
extended duals for SDPs and generalizations
 -- see Ramana \cite{Ramana:97}, and Klep and Schweighofer \cite{KlepSchw:12}. 
For the connection of these approaches, see 
Ramana, Tun\c{c}el and Wolkowicz \cite{RaTuWo:97}, and \cite{Pataki:13}. 
The third approach is that of elementary reformulations of SDPs, which is more recent -- see  \cite{Pataki:17} 
and  \cite{LiuPataki:15}. 

The reason that 
\eref{d}, $(\palt)$  and $(\dalt)$ are not exact  is that  
the linear image of a closed convex cone is not always closed.
For studies on when this image is closed (or not), see Bauschke and Borwein 
\cite{BausBor:99}; and Pataki \cite{Pataki:07}.

\co{
Borwein and Moors \cite{BorweinMoors:09, BorweinMoors:10} for a proof that 
``most'' linear transformations (in a well-defined sense) 
map a closed convex cone to a closed set. }

Here we unify, simplify, and extend the above approaches and develop a collection of 
exact duals, and certificates of infeasibility and weak infeasibility in conic LPs with the following features:

\benum
\item They do not rely on any constraint qualification (CQ), such as strict feasibility of \eref{p} (which requires that there exist $x \in \rad{m}$ such that 
$b - A x$ is in the relative interior of $K$).

\item They inherit most of the simplicity 
of the Lagrange dual (see Sections \ref{sec-fr},  \ref{sec-catalog}, and \ref{sec-catalog-sdp}).
Some of our infeasibility certificates generalize the row echelon 
form of a linear system of equations, as they consist of a small,
trivially infeasible subsystem obtained by elementary row operations. 
The size of this subsystem is bounded by a geometric parameter of the cone, the length of the 
longest chain of nonempty faces. 
The results for weak infeasibility are analogous.

\co{The ``easy'' proofs, as weak duality, and
the proofs of infeasibility and weak infeasibility are  nearly as simple as proofs 
in linear programming duality (see Sections \ref{sec-fr} and \ref{sec-catalog}). } 

\item Our exact dual of 
\eref{p} generalizes the exact SDP duals of Ramana \cite{Ramana:97} 
and Klep and Schweighofer \cite{KlepSchw:12} to the context of general conic linear programming:
it is a  conic LP whose constraints
are copies of the original constraints.

Our exact certificate of infeasibility of \eref{p} refines the certificate given by Waki and Muramatsu in \cite{WakiMura:12} by showing that it can be viewed as 
a  conic linear system.

\co{, and the underlying 
the {\em facial reduction cone} -- a convex cone we introduce, and use 
to encode facial reduction algorithms. }

\item They yield some fundamental geometric results in convex analysis: bounds on the number of constraints that can be 
dropped or added in a conic linear system while keeping it (weakly) infeasible 
(Corollary \ref{dimension} in Section \ref{sec-catalog}); 
an exact characterization of when the linear image of a closed convex cone is closed;  
and an exact characterization of  an important class of cones, called {\em nice cones}  
(see Section \ref{sec-geometry}).

\item They provide algorithms to generate {\em all} infeasible conic LP instances 
over several important cones (see Sections \ref{sec-catalog}, \ref{sec-generate},
 and \ref{sect-conclude}), 
and {\em all} weakly infeasible SDPs in a natural class (Section \ref{sec-generate}).

\item The above algorithms are easy to implement, and they provide a challenging test set 
of infeasible and weakly infeasible SDPs: while we can verify  the status of our instances by inspection in 
exact arithmetic, they are difficult for commercial and research codes (Section \ref{sec-comp}).

\item Of possible independent interest are  an 
elementary  facial reduction algorithm (Section \ref{sec-fr}) with a simplified 
proof of convergence  and the geometry of the {\em facial reduction cone}, a convex 
cone that we introduce  and use to encode facial reduction algorithms (see Lemma \ref{frkprop}). 

\eenum

We now describe our main tools, and some of our main results with full proofs of 
the ``easy'' directions. 
We will often reformulate a conic LP in a suitable form from which its status (as infeasibility) 
is easy to read off. This process is akin to bringing a matrix to row echelon form, and 
most of the operations we use indeed come from Gaussian elimination. 
To begin, we represent $A$ and $A^*$ as 
$$
Ax \, = \, \sum_{i=1}^m x_i a_i, \, A^*y \, = \, ( \la a_1, y \ra, \dots, \la a_m, y \ra)^T, \, \mbox{where $a_i \in Y \,$ for $i=1, \dots, m.$} 
$$
\co{
To begin, we will frequently write the primal-dual pair in the form 
\begin{center}                                                                                              
$$                                                                                                          
\begin{array}{lrlcrlr}                                                                                      
    &   \sup  & \sum_{i=1}^m c_i x_i          &    \hspace{2cm} & \inf   &   \la b, y \ra  &   \\ 
(P) &   s.t.  & \sum_{i=1}^m x_i a_i \leq_K b &    \hspace{2cm} & s.t.  &   \la a_i, y \ra = c_i \,\, \forall i & (D)     \\ 
    &         &                               &    \hspace{2cm} &       &   y \geq_{K^*} 0  &    \end{array}                       $$   
\end{center}                                                                                                      
where $a_i \in Y$ for all $i.$ }
\co{
\begin{Definition} \label{reform-def}
We obtain an {\em elementary reformulation} or {\em reformulation}
of $(P)\mhyphen(D)$ by a sequence of the operations: 1) multiply a dual equation $\la a_i, y \ra = c_i\,$ 
by a nonzero scalar; 2) add a multiple of 
a dual equation to another; 3) exchange two equations; 4) replace $b$ by $b + A \mu$ for some $\mu \in \rad{m}.$

If $K = K^*$ we also allow the operation: 5) replace $a_i$ by $T a_i  (i=1, \dots, m)\,$ and $b$ by $Tb, \,$ 
where $T$ is an invertible linear map with $TK = K.$ 
\end{Definition}
}
\co{
To begin, we represent $A$ and $A^*$ as 
$$
Ax \, = \, \sum_{i=1}^m x_i a_i, \, A^*y \, = \, ( \la a_1, y \ra, \dots, \la a_m, y \ra)^T, \, \mbox{where $a_i \in Y \,$ for $i=1, \dots, m.$} 
$$
}
\begin{Definition} \label{reform-def}
We obtain an {\em elementary reformulation} or {\em reformulation}
of $(P)\mhyphen(D)$ by a sequence of the operations:
\benum
\item \label{lambda} Replace $(a_{i}, c_{i})$ by $(A \lambda,\la c,\lambda\ra)$ for some $i \in \{1, \dots, m \},$
where $\lambda \in \rad{m}, \, \lambda_i \neq 0.$
\item \label{switch} Switch $(a_{i},c_{i})$ with $(a_{j},c_{j})$, where $i \neq j$ 
\footnote{Performing a sequence of operations of type 
	(\ref{lambda}) and (\ref{switch}) 
	is the same as replacing $A$ by $AM$ and $c$ by $M^Tc,$ where $M$ is an $m \ti m$ invertible matrix.}.  
\item \label{mu} Replace $b$ by $b+A\mu$, where $\mu \in \rad{m}$.
\eenum
If $K = K^*$ we also allow the operation:
\benum
\item[(4)] Replace $a_i$ by $T a_i  (i=1, \dots, m)\,$ and $b$ by $Tb, \,$ 
where $T$ is an invertible linear map with $TK = K.$ 
\eenum
\end{Definition}
We call operations (1)-(3) {\em elementary row operations}. 
Sometimes we reformulate only \eref{p} or \eref{d}, or only the underlying 
systems, ignoring the objective function. 
Clearly, a conic linear system is infeasible, strongly infeasible, etc., exactly when 
its elementary reformulations are. 

{\em Facial reduction cones} ``encode'' a facial reduction 
algorithm, in a sense that we make precise later, and 
will replace the usual dual cone to make our duals, and infeasibility certificates exact.  
\begin{Definition} \label{frk-def} 
Let $k \geq 0$ be an integer. 
The order $k$ facial reduction cone of $K$ is the set 
\beqast
\fr_k(K) = \{\, (y_{1},\ldots,y_{k}) \, : \, y_1 \in K^*, \, y_{i}\in(K\cap y_{1}^{\perp}\cap \ldots\cap y_{i-1}^{\perp})^{*}, i =2,\ldots k \, \}.
\eeqast
We drop the subscript  $k$ when its value is clear from the context  or if it is irrelevant. 
\end{Definition}
We have $\fr_0(K) = \emp, \, \fr_1(K) = K^*$ and we can pad elements of $\fr_k(K)$ with zeros to make them elements 
of a higher order facial reduction cone,  hence these  sets 
serve as  relaxations of $K^*.$ 
Surprisingly, $\fr_k(K)$ is convex, which is closed only in trivial cases, but 
behaves as well as 
$K^*$ under the usual operations on convex sets -- see Lemma \ref{frkprop}. 

We now state an excerpt of our main results with full proofs of the ''easy'' directions: 

\nin{\bf Theorem I} If $K$ is a general closed convex cone, then 
\benum 
\item \label{thmI1} \eref{d} is infeasible if and only if it has a reformulation 
\beq\label{dref}\tag{$D_{\rm ref}$}
\ba{rcl}
\la a_{i}^{\prime}, y\ra & = & 0\,(i=1,\ldots,k)\\
\la a_{k+1}^{\prime}, y\ra & = & -1\\
\la a_{i}^{\prime}, y\ra & = & c_{i}^{\prime}\,(i=k+2,\ldots,m)\\
y & \geq_{K^*} & 0
\ena
\eeq
where $k \geq 0, \, (a_{1}^{\prime},\ldots,a_{k+1}^{\prime})\in\fr(K^*)$.
\item \label{thmI2} \eref{d} is not strongly infeasible if and only if there is $\ell \geq 0$ and 
$(y_1, \dots, y_{\ell+1}) \in \fr(K)$, such that
\beqast
\ba{ccl}
A^*y_{i} &=& 0 \,(i=1,\ldots,\ell)\\
A^*y_{\ell+1} &=& c.
\ena
\eeqast
\eenum
\qed

Part (1) of Theorem I generalizes the row echelon form of a linear system of 
equations: it finds a small, trivially infeasible subsystem in \eref{d} using only elementary row
operations. As we prove later, the size of the subsystem (i.e., $k+1$) 
is
bounded by a geometric parameter of $K^*$. Also, if  $K^*$ is the whole space, then 
$\fr_{k+1}(K^*) = \{0\}^{k+1}, \,$ so the constraint $\la 0, y \ra = -1 $ in \eref{dref}
proves infeasibility. 
\co{When \eref{d} is strongly infeasible, 
	we can take $k = 0$ in part \eref{thmI1} with $a_1^\prime \in K, \,$ and 
	$(a_1^\prime, -1)  = (Ax, \la c, x \ra)$ for some  $x. \,$ 
	In the general case 
	$k$ will be bounded by the length of the longest chain of faces in $K^*$ (see Section \ref{sec-catalog}). }
Turning now to part (2), 
we note that	if $\ell = 0$ then \eref{d} is actually feasible. 

\co{Turning now to part (2) we note that
	if $\ell = 0$ then \eref{d} is actually feasible. }

Naturally, 
combining the two parts of Theorem I we obtain an exact characterization of weak infeasibility.
We prove the ''easy'', i.e., the ''if'' directions below: 

\pf{of ''if'' in part (1)}
We will prove that \eref{dref} is infeasible, so suppose that $y$ is feasible in it to obtain the 
contradiction 
$$
y\in K^*\cap a_{1}^{\prime\perp}\cap\ldots\cap a_{k}^{\prime\perp} \Rightarrow \la a_{k+1}^{\prime}, y \ra \geq 0.
$$

\pf{of ''if'' in part (2)} Let us fix $(y_1, \dots, y_{\ell+1})$ as stated, and 
$x$ such that  $Ax \geq_K 0.$ 
Then 
$$
\la c, x \ra = \la A^*y_{\ell+1}, x \ra= \la  y_{\ell+1}, Ax \ra \geq 0,
$$
where in the inequality we used $Ax \in K \cap y_1^\perp \cap \dots \cap y_\ell^\perp.$ Thus 
$(\dalt)$  cannot be feasible. \qed

We illustrate  Theorem I with a semidefinite system, with 
$Y = \sym{n} \,$ the set of order $n$ symmetric matrices and  
$K = K^* = \psd{n}$ as the set of 
order $n \,$ symmetric positive semidefinite  (psd) matrices.
The inner product of $a, b \in \sym{n}$ is
$a \bullet b := \la a, b \ra := \trace(ab) \,$ 
and we write $\preceq$ in place of $\leq_K.$ 
Note that we denote the elements of $\sym{n}$ by small letters, and we reserve capital letters 
for operators.

\bex \label{ex1} {\rm 
The semidefinite system 
\beq \label{ex1-sys}
\ba{rcl} 
\bpx 1 & 0 \\ 0 & 0 \epx \bullet y & = & 0 \\
\bpx 0 & 1 \\ 1 & \alpha \epx \bullet y & = & -1 \\
                                      y & \succeq & 0    
\ena
\eeq
is infeasible for any $\alpha \geq 0,$ and weakly infeasible exactly when $\alpha =0.$ 

Since the constraint matrices in \eref{ex1-sys} are in $\fr(\psd{2}),$ 
this system is in the form of \eref{dref} (and itself is a proof of infeasibility).

Suppose $\alpha=0$ and let 
\beq \label{y1y2} 
y_1 = \bpx 0 & 0 \\ 0 & 1 \epx, \, y_2 = \bpx 0 & -1/2 \\ -1/2 & 0 \epx.
\eeq
Then $(y_1, y_2) \in \fr(\psd{2}), \,A^*y_1 = (0,0)^T, \, A^*y_2 = (0, -1)^T, \,$ hence 
$(y_1, y_2)$ proves that  \eref{ex1-sys} is not strongly infeasible. 

For a set $S$ we denote by $\front(S)$ its {\em frontier}, i.e., the difference between its closure, and itself:
\beq
\front(S) := \cl S \setminus S.
\eeq
Classically,  the set of right hand sides that make \eref{d} weakly infeasible is the 
frontier of $A^*K^*;$  in the $\alpha = 0$ case 
\beq \label{front-eq} 
\front(A^*K^*) \, = \, \, \{ \, (0, \lambda)^T: \, \lambda \neq 0 \, \} \footnote{To outline a proof of 
	(\ref{front-eq}), assume for simplicity $\lambda=2. \,$ Then 
	$$
	y := \bpx \eps & 1  \\
	1 & 1/\eps \epx \succeq 0 \; {\rm for \, all} \,  \eps > 0, 
	$$
	and $A^*y = (\eps, 2)^T, \,$ 
	so  $(0,2)^T \in \cl (A^* \psd{2});$ the rest  of (\ref{front-eq}) is straightforward to verify.}; 
\eeq
For all such right hand sides a suitable $(y_1,y_2) \in \fr(\psd{2})$ 
proves that \eref{ex1-sys} is not strongly infeasible.
}
\eex

We organize the rest of the paper as follows. In 
the rest of the introduction we review prior work, collect notation, and record basic properties of the facial reduction cone $\fr_k(K)$. In Section  \ref{sec-fr} we present our simple facial 
reduction algorithm, and our exact duals of \eref{p} and \eref{d}. 
In Section \ref{sec-catalog} we present our exact certificates of infeasibility and 
weak infeasibility of general conic LPs. 
 In Section \ref{sec-catalog-sdp} we describe 
specializations to SDPs. 
Section \ref{sec-geometry} presents our  geometric corollaries: an exact characterization of 
when the linear image of a closed convex cone is closed, and an exact characterization 
of nice cones. 
In Section \ref{sec-generate} we give  our algorithm to generate {\em all}
infeasible SDPs, define a natural class of weakly infeasible SDPs, and 
provide a simple algorithm to generate {\em all} instances in this class. In Section 
\ref{sec-comp} we present a library of infeasible, and weakly infeasible SDPs, and 
our computational results. In Section \ref{sect-conclude} we discuss possible extensions of our work,
and conclude the paper. 

The main ideas of the paper can be  quickly absorbed by reading only Sections \ref{sec-fr}, \ref{sec-catalog} and 
\ref{sec-catalog-sdp}. A reader interested in the geometrical/convex analysis aspects will probably want to read 
Section \ref{sec-geometry}; other readers can skip this section at first reading.
Section \ref{sec-generate} gives more insight into the structure of weakly infeasible 
SDPs, and Section \ref{sec-comp}
is more for an audience interested in computation.  
The paper relies only on knowledge of  elementary convex analysis and the results are illustrated by many examples. 

\noindent{\bf Prior work} 
Facial reduction algorithms turn \eref{d} into an exact dual of 
\eref{p} by constructing a suitable smaller cone  $F$, replacing $K$ by $F$ in \eref{p}, and 
$K^*$ by $F^*$ (a larger cone) in \eref{d}. 
The first such algorithm was proposed by  Borwein and Wolkowicz \cite{BorWolk:81,BorWolk:81B} for 
nonlinear conic systems. Their algorithm assumes that the system is feasible. 
Waki and Muramatsu \cite{WakiMura:12} described a variant for conic linear systems, 
which is simpler, and it is the first variant of facial reduction, which 
allows one to prove infeasibility. Pataki \cite{Pataki:13} proposed another simplified version, which is 
fairly straightforward to implement. 

We can construct exact duals of SDPs without relying on any constraint qualification, at the cost of
introducing
polynomially many extra variables and constraints --  see Ramana \cite{Ramana:97}, 
and Klep and Schweighofer \cite{KlepSchw:12}. 
 We note that Ramana's dual 
relies on convex analysis, while  Klep and Schweighofer's dual uses ideas from algebraic geometry, namely sums of  
squares representations. 
While at first they seem different, 
the approaches of facial reduction and  extended duals 
are related -- see Ramana, Tun\c{c}el and Wolkowicz \cite{RaTuWo:97}, and \cite{Pataki:13} for proofs
of the correctness of Ramana's dual relying on facial reduction algorithms.
The paper \cite{Pataki:13} generalizes Ramana's dual to the context of conic LPs over {\em nice} cones,
while P\'{o}lik and Terlaky \cite{PolikTerlaky:09} extend Ramana's dual to conic LPs 
over homogeneous cones. Ramana and Freund in \cite{RaFreund:96} studied the Lagrange dual of the extended dual of Ramana in \cite{Ramana:97}, and proved that it has the same optimal value as the original
problem.  

For recent studies on the closedness of the linear image of a closed convex cone
we refer to Bauschke and Borwein \cite{BausBor:99};   and Pataki \cite{Pataki:07}.
The paper \cite{BausBor:99} gives a necessary and sufficient condition for the continuous image of a closed convex cone 
to be closed, in terms of the {\em strong conical hull intersection property.} 
Pataki \cite{Pataki:07} gives necessary conditions, which subsume well known sufficient conditions, and 
are necessary and sufficient for a broad class of cones, called nice cones. 
See also \cite[Lemma 1]{Pataki:17}) for a simplified exposition. 
We refer to Bertsekas and Tseng
\cite{BertTseng:07} for a study of a more general problem, 
whether the intersection of a nested sequence of closed sets is nonempty. 
The paper \cite{Pataki:17} applies the closedness result of \cite{Pataki:07} to 
characterize when a conic linear system is {\em badly behaved;} by bad behavior of  $Ax \leq_K b$ we mean that there is a duality gap between (\ref{p}) and 
(\ref{d}),  or (\ref{d}) does not attain its optimal value {\em for some} $c.$

Borwein and Moors \cite{BorweinMoors:09, BorweinMoors:10} recently showed that the 
set of linear maps under which the image is {\em not} closed is $\sigma$-porous, i.e., it 
has Lebesgue measure zero, and is also small in terms of category.
For characterizations of nice cones, see Pataki \cite{Pataki:12}; 
and Roshchina \cite{Vera:13} for a proof that not all facially exposed cones are nice. 

Elementary reformulations of  SDPs -- see Pataki \cite{Pataki:17} 
and Liu and Pataki \cite{LiuPataki:15} -- use simple operations, as 
elementary row operations, to bring a semidefinite 
system into a form from which its status (as infeasibility) is trivial to read off. 
Definition \ref{reform-def} generalizes  elementary reformulations of SDPs
to the context of general conic LPs.  

\co{The approaches of facial reduction, extended duals and elementary reformulations are all related -- see 
Ramana, Tun\c{c}el and Wolkowicz \cite{RaTuWo:97}, \cite{Pataki:13} 
and Liu and Pataki \cite{LiuPataki:15}.}
Another related paper is by Lourenco et al \cite{Lourenco:13} which presents 
an error bound based reduction procedure  to simplify weakly infeasible SDPs, and a proof 
that all such SDPs with order $n$ constraint matrices 
contain a weakly infeasible subsystem with dimension at most $n-1.$  
Part (4) of our Theorem \ref{mainfarkas} generalizes this result to the context of general conic LPs. 

For an application of facial reduction to the {\em euclidean distance matrix completion problem} see 
Krislock and Wolkowicz \cite{KrisWolk:10}; and Drusviyatsky et al \cite{Dima:15} for 
a theoretical analysis  of the algorithm in \cite{KrisWolk:10}. 
We refer to  
Waki in \cite{Waki:12} for a method to generate weakly infeasible SDPs from Lasserre's 
relaxation of polynomial optimization problems.

For textbook treatments of the duality theory of conic LPs we refer to 
Bonnans and Shapiro \cite{BonnShap:00}; Renegar \cite{Ren:01};  G$\ddot{\mathrm{u}}$ler \cite{Guler:10}; and Borwein and Lewis \cite{BorLewis:00}. 

A remark on notation: we consider the primal problem to 
be in inequality form, since this is how Ramana \cite{Ramana:97}, and Klep and Schweighofer 
\cite{KlepSchw:12} presented their
dual, and infeasibility 
certificate. The equality constrained dual form will
be more useful to derive our geometric corollaries (in Section \ref{sec-geometry}) and to generate 
instances (in Section \ref{sec-generate}). 

\noindent{\bf Notation and preliminaries}
We assume throughout that the operator $A$ is surjective. 
For an operator $M$ we denote by $\R(M)$ its  rangespace, and by 
$\N(M)$ its nullspace.
For $x$ and $y$ in the  same Euclidean space we sometimes 
write $x^*y$ for $\la x, y \ra.$ 
For a convex set $C$ we denote its linear span, the orthogonal complement of 
its linear span, its closure, and relative interior by 
$\lin C, \, C^\perp, \, \cl C, \,$ and $\ri C, \,$ respectively.
We define the dual cone of $K$ as 
$$
K^* \, = \, \{ \, y \, | \, \la y, x \ra \geq 0, \, \forall x \in K \, \},
$$
and for convenience we set 
$$ \label{defperp}
K^{* \setminus \perp} \, := \, K^* \setminus K^\perp.
$$
\co{We say that \eref{p} is an exact dual of \eref{d} if it has the same value, as 
\eref{d}, and attains this value, when it is finite. 
We know that \eref{p} is an exact dual of \eref{d} if the latter is strictly feasible, i.e., 
if there is $\by \in \ri K^*$ feasible in \eref{d}.}
 
 We say that \eref{d} is {\em strictly feasible}, 
 if there is $y \in \ri K^*$ feasible in it. If \eref{d} is strictly feasible, 
 then \eref{p} is an exact dual of \eref{d} (i.e., it has the same value, as 
 	\eref{d}, and attains this value, when it is finite).

For $F, \,$  a convex subset of $K$ we say that 
$F$ is a {\em face} of $K, \,$ 
if $y, z \in K, \,$ and $1/2(y+z) \in F$ implies $y, z \in F.$ 

We denote by $\ell_K$ the length of the longest chain of nonempty faces in $K, \,$ i.e.,
\beq \label{defellK}
\ell_K \,  = \, \max \, \{ \, k \, | \, F_1 \subsetneq F_2 \subsetneq \dots \subsetneq F_k \, {\rm are \, nonempty \,  faces \, of \,} K \, \}.
\eeq
For instance, $\ell_{\psd{n}} = \ell_{\rad{n}_+} = n+1. \,$ 

If $H$ is an affine subspace with $H \cap K \neq \emptyset, \,$ 
then we call the smallest face of $K$ that contains $H \cap K$  the {\em minimal cone of} 
$H \cap K, \,$ i.e., it is
\beq \label{mincone-eq} 
\bigcap \{ \, G \, | \, G \supseteq H \cap K, \, G \, {\rm is \, a \, face \, of \,} K \}. 
\eeq

We call $z \in K$ a {\em slack} in \eref{p}, if 
$z = b - Ax$ for some $x \in \rad{m};$ if 
\eref{p} is an SDP, then we call such a $z$ a {\em slack matrix.} 

\co{
The traditional alternative system of $(P)$ is 
\beq \label{palt}\tag{$P_{\rm alt}$}
\ba{rlcl}
A^{*}y          &=   & 0  \\
b^* y           & =  & -1       \\
 y              & \geq_{K^*}  &  0;
\ena
\eeq
if   it is feasible, then \eref{p} is infeasible and we say that it is strongly infeasible.
If \eref{p} and \eref{palt} are both infeasible, then we say that \eref{p} is weakly infeasible.  
}

For a nonnegative integer $r$ we denote by $\psd{r} \oplus \{ 0 \}$ the subset of 
$\psd{n}$ (where $n$ will be clear from the context) with psd upper left $r \ti r$ block, and the rest zero.
All faces 
of $\psd{n}$ are of the form $t^T (\psd{r} \oplus 0)t$ where $t$ is an invertible matrix 
\cite{BarCar:75, Pataki:00A}. 
We write 
\beqa \label{def-autk}
\aut(K) & = & \{ \, T: Y \rightarrow Y \, | \, T \, {\rm is \, linear \, and \, invertible,  \,} T(K) = K \, \}
\eeqa
for the automorphism group of a closed convex cone $K.$ 
For example,  $T \in \aut(\psd{n})$ if and only if there is 
an invertible matrix $t$ such that 
$T(x) = t^T x t$ for $x \in \sym{n}.$ 

\co{
In this section we record some elementary properties of $\fr_k(K).$ Somewhat surprisingly,
$\fr_k(K)$ is a convex cone, which is only closed in trivial cases. Despite this shortcoming it behaves 
nearly as well as the usual dual cone $K^*$ under the usual operations on convex sets. For instance
any linear transformation that preserves $K$ also 
preserves $\fr_k(K).$ We note that part \eref{frkprop3} is stated with some abuse of notation 
(see the proof). 
}

In Lemma  \ref{frkprop} we record relevant properties of $\fr_k(K).$ 
Its proof is given in Appendix \ref{app-proofs}.
\ble \label{frkprop}
For $k \geq 1 \,$ the following hold: 
\benum
\item \label{frkprop1} $\fr_k(K)$ is a convex cone.
\item \label{frkprop15} $\fr_k(K)$  closed if and only if $K$ is a subspace or $k=1.$ 
\item \label{frkprop2} If $T \in \aut(K),$ and $(y_1, \dots, y_k) \in \fr_k(K) \,$ then 
$$
(T^*y_1, \dots, T^*y_k) \in \fr_k(K).  
$$
\item \label{frkprop3} If $C$ is another closed convex cone, then 
\beq \label{abuse} 
\fr_k(K \times C) \, = \, \fr_k(K) \times \fr_k(C).
\eeq
Precisely, 
\beq \label{KC}
\bigl( (y_1, z_1), \dots, (y_{k}, z_{k})\bigr) \, \in \, FR_k(K \times C)
\eeq
if and only if 
\beq \label{KC2}
\bigl(y_1, \dots, y_{k}\bigr) \in \fr_k(K) \, {\rm and} \, \bigl(z_1, \dots, z_{k}\bigr) \in \fr_k(C).
\eeq
\eenum
\ele
Note that equation  (\ref{abuse}) somewhat abuses notation; the statement 
$(\ref{KC}) \LRA (\ref{KC2})$ is fully rigorous. 
	Also note that with $k=1$  statement (\ref{frkprop2}) recovers the well known fact:
	$T \in \aut(K) \Rightarrow T^* K^* \subseteq K^*.$

\section{Facial reduction and exact duals in conic linear programming}
\label{sec-fr}

In this section we present a very simple facial reduction algorithm to 
find $F, \, $ the minimal cone of the system 
\beq \label{HK} 
H \cap K,
\eeq
where $H$ is an affine subspace  with $H\cap K\neq \emptyset, \,$ and our 
exact duals of \eref{p} and of \eref{d}. 
The convergence proof of Algorithm \ref{algo:FRA}, with an upper bound 
on the number of necessary steps, is entirely elementary, and it simplifies the 
proofs given in \cite{WakiMura:12} and \cite{Pataki:13}.  Recall the definition of the 
minimal cone from \eref{mincone-eq}.) 

We first illustrate the importance of the minimal cone: 
if $F$ is the minimal cone of \eref{HK} then $$ H \cap \ri F \neq \emptyset$$ (otherwise $H \cap K$ would be 
contained in a proper face of $F$). Assume next 
that $F$ is the minimal cone of $(\R(A) + b) \cap K,$ where $A$ and $b$ are as in \eref{p}. Then 
replacing $K$ by $F$ makes \eref{p} strictly feasible, and keeps its feasible set the same.
Hence if we also replace $K^*$ by $F^*$ in \eref{d},  then \eref{d} becomes 
an exact dual of \eref{p}. 

As illustration, we consider the following example:
\bex \label{a1a2a3-ex} {\rm The optimal value of the SDP 
\beq \label{a1a2a3-sys}
\ba{rl}
\sup & x_1 \\
s.t. & x_1 \bpx 0 & 1 & 0 \\ 1 & 0 & 0 \\ 0 & 0 & 0 \epx + x_2 \bpx 0 & 0 & 1 \\ 0 & 1 & 0 \\ 1 & 0 & 0 \epx \preceq \bpx 1 & 0 & 0 \\ 0 & 0 & 0 \\ 0 & 0 & 0 \epx
\ena
\eeq
is zero. Its usual SDP dual, in which we denote the dual matrix 
by $y$ and its components by $y_{ij}$, is equivalent to 
\beq \label{supx1-d} 
\ba{rl}
\inf & y_{11} \\
s.t. & \bpx y_{11} & 1/2        & - y_{22}/2 \\
             1/2     & y_{22}   & y_{23} \\
        - y_{22}/2     &  y_{23}  & y_{33} \epx \succeq 0,
\ena
\eeq
which does not have a feasible solution with $y_{11} = 0$ (in fact it has an unattained $0$ infimum). 

Since all slack matrices in \eref{a1a2a3-sys} 
are contained in $\psd{1} \oplus \{ 0 \}$ and there is a slack matrix whose 
$(1,1)$ element is positive, the  minimal cone of this system is 
$F \, = \, \psd{1} \oplus \{ 0 \}.$ 
If in the dual program we replace $\psd{3}$ by $F^*,$ then 
the new dual  attains with 
\beq \label{yf*}
y \, := \, \bpx 0   & 1/2 & 0 \\
           1/2 & 0   & 0 \\
            0  & 0   & 0 
       \epx \, \in \, F^* \, = \, \bpx \oplus & \ti & \ti \\
                        \ti    & \ti & \ti \\
                        \ti    & \ti & \ti \epx \, 
\eeq
being an optimal solution. Here the $\oplus$ and $\ti$ symbols stand for  psd and arbitrary submatrices, respectively.
}
\eex
To construct the minimal cone of $H \cap K$ 
we rely on the following classic theorem of the alternative 
(recall $K^{* \setminus \perp} = K^* \setminus K^\perp$). 
\beq \label{HKlemma}
H\cap \ri{K}=\emptyset\Leftrightarrow H^{\perp}\cap K^{*\setminus \perp}\neq \emptyset.
\eeq
Algorithm \ref{algo:FRA} repeatedly applies the equivalence \eref{HKlemma} to find $F:$ 
\begin{algorithm}[H] 
\caption{Facial Reduction}
\label{algo:FRA}
\begin{algorithmic}
\INITIALIZATION Let $y_0 = 0, \, F_0 = K, \, i=1.$
\While {$\exists y_i \in H^\perp \cap F_{i-1}^{* \setminus \perp}$}
\State Choose such a $y_{i}.$
\State Let $F_i = F_{i-1} \cap y_i^\perp.$
\State  Let  $i = i+1.$
\EndWhile
\end{algorithmic}
\end{algorithm}
Algorithm \ref{algo:FRA} generates a sequence of faces
\beq \label{FFk}
F = F_k \subsetneq F_{k-1} \subsetneq \dots F_0 = K.
\eeq
The reason that $F = F_k$ holds for some $k$, is that once we stop, 
$\ri F_k \cap H \neq \emp. \,$ Since  
$$
\ri F_k \cap H \, = \, \ri F_k \cap H \cap K \, = \, \ri F_k \cap H \cap F,  
$$
we deduce $\ri F_k \cap F \neq \emp.$ Thus by Theorem 18.1 in \cite{Rockafellar:70} 
we obtain $F_k \subseteq F $ with the reverse containment already given. 
\co{The reason that $F = F_k$ holds for some $k$, is that once we stop, 
$\ri F_k \cap H \neq \emp. \,$ Hence 
$$
\ri F_k \cap H \cap K \, = \, \ri F_k \cap H \cap F \, \neq \, \emp,
$$
so $\ri F_k \cap F \neq \emp.$ Thus by Theorem 18.1 in \cite{Rockafellar:70} 
we obtain $F_k \subseteq F $ with the reverse containment already given. } 
\co{
The reason that $F = F_k$ holds for some $k$, is that once 
we stop, 
$$
\emptyset \neq \ri F_k \cap H \, = \, \ri F_k \cap H \cap K \, = \, \ri F_k \cap H \cap  F \subseteq \ri F_k \cap F,
$$
hence by Theorem 18.1 in \cite{Rockafellar:70} 
we obtain $F_k \subseteq F $ with the reverse containment already given. 
}

To bound the number of steps that Algorithm \ref{algo:FRA} needs,  we need a definition: 
\begin{Definition} For $k \geq 1$ we say that $(y_{1},\ldots,y_{k}) \in \fr_k(K)$ is {\em strict}, if 
$$
y_i \, \in \, (K \cap y_1^\perp \cap \dots \cap y_{i-1}^\perp)^{* \setminus \perp} \, {\rm for} \, i=2, \dots, k.
$$
We say that it is {\em pre-strict} if $(y_{1},\ldots,y_{k-1})$ is strict. 
\end{Definition}
\bth \label{frlemma}  If $(y_{1},\ldots,y_{k})$ is strict, then these vectors are linearly independent.
Also, Algorithm \ref{algo:FRA} stops after at most 
$$ 
\min \, \{ \, \ell_K - 1, \, \dim H^\perp \, \}
$$ 
iterations.
\enth
\pf{} To prove the first statement, assume the contrary. Then for some $2 \leq i \leq k$ a  contradiction follows:
$$
y_{i} \in \lin \, \{ \, y_1, \dots, y_{i-1} \, \} \, \subseteq \, (K \cap y_1^\perp \cap \dots \cap y_{i-1}^\perp)^\perp.
$$
The second statement is then immediate.
\qed

\bex{\rm (Example \ref{a1a2a3-ex} continued)
Algorithm 1 applied to this example may output the sequence 
\renewcommand{\arraystretch}{1.1}
\begin{equation} \label{y1y2}
\begin{array}{rcccl} 
y_1 & = & \begin{pmatrix} \,0 & \,\;0\,\; & 0\, \\ \,0 & \,\;0\,\; & 0\, \\ \,0 & \,\;0\,\; & 1\, \end{pmatrix}&,& \, F_1 \, = \, 
\begin{pmatrix} 
\parbox{1cm}{$\,\, \bigoplus$}  & \hspace{-.25cm} \parbox{1cm} {$\begin{array}{c} 0 \\ 0 \end{array}$} \hspace{-.55cm} \\
\parbox{1cm}{$0 \,\,\,\, 0 \,\,\,\,$}  &   \hspace{-.5cm} 0 
\end{pmatrix}, \,\, \\
y_2 & = & \begin{pmatrix} 0 & 0 & -1 \\ 0 & 2 & 0 \\ -1 & 0 & 0 \end{pmatrix}&,& \, F_2 \, = \, F, 
\end{array}
\end{equation}
\renewcommand{\arraystretch}{1}
\co{
\renewcommand{\arraystretch}{1.1}
\begin{equation} \label{y1y2}
\begin{array}{rcl} 
y_1 & = & \begin{pmatrix} 0 & 0 & 0 \\ 0 & 0 & 0 \\ 0 & 0 & 1 \end{pmatrix}, \, F_1 \, = \, 
\begin{pmatrix} 
  \parbox{1cm}{$\,\, \bigoplus$}  & \hspace{-.5cm} \parbox{1cm} {$\begin{array}{c} 0 \\ 0 \end{array}$} \hspace{-.5cm} \\
  \parbox{1cm}{$0 \,\,\,\, 0 \,\,\,\,$}  &   \hspace{-.45cm} 0 
\end{pmatrix}, \,\, \\
y_2 & = & \begin{pmatrix} 0 & 0 & -1 \\ 0 & 2 & 0 \\ -1 & 0 & 0 \end{pmatrix}, \, F_2 \, = \, F,
\end{array}
\end{equation}
\renewcommand{\arraystretch}{1}} 
where again the $\oplus$ symbol corresponds to a psd submatrix.
Note that $y_i \in H^\perp$ in this context means 
that $y_1$ and $y_2$ are orthogonal to all constraint matrices in the system (\ref{a1a2a3-sys}).  
}
\eex
From Theorem \ref{frlemma} we immediately obtain an extended exact
 dual for \eref{p}, described below in 
\eref{Dstrong}. Note that 
\eref{Dstrong} is a  conic linear program whose data is the same as the data of \eref{p}, thus it extends the 
exact SDP duals of Ramana \cite{Ramana:97} and Klep and Schweighofer \cite{KlepSchw:12} to the context 
of general conic LPs. The underlying cone in \eref{Dstrong} is the facial reduction cone: 
thus, somewhat counterintuitively, we find an exact dual of \eref{p} over a 
convex cone, which is not closed in all important cases (see Lemma \ref{frkprop}). 

\bth \label{sdlemma} For all large enough $k$ the problem 
\beq\label{Dstrong}\tag{$D_{\rm ext}$}
\ba{rccl}
\inf & b^* y_{k+1} &&\\
s.t. &A^*y_{k+1}&=&c\\ 
     & A^{*}y_{i} &=& 0\, (i=1,\ldots,k) \\ 
     &  b^* y_i  & = & 0 \, (i=1,\ldots,k) \\
     &(y_{1},\ldots,y_{k+1})&\in& \fr_{k+1}(K)
\ena
\eeq
is an exact dual of $(P).$ 

When the value of \eref{p} is finite, for a suitable $k$ the problem 
\eref{Dstrong} has a pre-strict optimal solution with the same value. 

\co{When the value of \eref{p} is finite, \eref{Dstrong} has an optimal solution 
$(y_1, \dots, y_{k'}, 0, \dots, 0, y_{k+1} pre-strict optimal solution with the same value. }

\enth 
\pf{} We first prove weak duality. 
Suppose that $x$ is feasible in $(P)$ and $(y_1, \dots, y_{k+1})$ in \eref{Dstrong},  then
\beqast
\la b, y_{k+1} \ra - \la c, x \ra & = & \la b, y_{k+1} \ra - \la A^* y_{k+1}, x \ra \\
                                  & = & \la b - Ax, y_{k+1} \ra \geq 0, 
\eeqast
where the last inequality follows from $b - Ax \in K \cap y_1^\perp \cap \dots \cap y_k^\perp.$ 

To prove the rest of the statements, first assume 
that \eref{p} is unbounded. Then by weak duality \eref{Dstrong} is
infeasible. Suppose next that \eref{p}  has a finite value $v, \, $ and 
let $F$ be the minimal cone of \eref{p}.  
Let us choose $y \in F^*$ to satisfy the affine constraints of \eref{d} with $b^*y=v.$ 
We have that 
$$
F   \, = \, K \cap y_1^\perp \cap \dots \cap y_k^\perp \, 
$$
for some $k \geq 0$ and $(y_1, \dots, y_k) \in \fr_k(K)$ strict, with all $y_i$ in 
$(\R(A)+b)^\perp.$ Hence if we choose this particular $k$ in \eref{Dstrong}, then 
$(y_1, \dots, y_k, y)$ is feasible in it with value $v.$ (If we choose $k$ larger, we can just pad
the sequence of $y_i$ with zeros.) This completes the proof. 
\qed

\bex \label{continued} {\rm (Example \ref{a1a2a3-ex} continued)
If we choose $y_1, y_2$ 
as in \eref{y1y2} and $y$ as in \eref{yf*}, then 
$(y_1, y_2, y)$ is an optimal solution to the extended dual of 
\eref{a1a2a3-sys}. 
}
\eex

\bth \label{splemma}
If \eref{d} is feasible then it has a strictly feasible reformulation 
\beq\label{dreffeas} \tag{$D_{\rm ref, feas}$}
\ba{rrcl}
\inf & b^* y \\
s.t. & \la a_{i}^{\prime}, y\ra & = & 0\,(i=1,\ldots,k)\\
& \la a_{i}^{\prime}, y\ra & = & c_{i}^{\prime}\,(i=k+1,\ldots,m)\\
& y & \in & K^* \cap a_1^{\prime \perp} \cap \dots \cap a_k^{\prime \perp},
\ena
\eeq
with $k \geq 0, \, (a_{1}^{\prime},\ldots,a_{k}^{\prime})\in\fr(K^*), $ 
which can be chosen strict. 
\enth
\pf{} Let us fix $\by$ such that $A^* \by = c, \,$  
and let 
$G$ be the minimal cone of 
$$
(\by + \N(A^*)) \cap K^*,
$$
i.e., of the feasible set of \eref{d}. Since Algorithm 
1 can construct $G, \,$ there is 
$k \geq 0$ and a strict 
$\, (a_{1}^{\prime},\ldots,a_{k}^{\prime})\in\fr(K^*)$ such that
$$
\ba{rcll}
a_{i}^\prime   & \in & \R(A) \cap \by^\perp \,  (i=1, \dots, k), \\
G              &  =  & K^* \cap a_{1}^{\prime \perp} \cap \dots \cap a_{k}^{\prime \perp}.
\ena
$$
By Theorem \ref{frlemma}  the vectors $a_1^\prime, \dots, a_k^\prime$ are linearly independent, so 
we can expand them to 
$$
A^\prime = [a_1^\prime, \dots, a_k^\prime, a_{k+1}^\prime, \dots, a_m^\prime], \, {\rm a \, basis \, of \, } \R(A).
$$
Let us write $A^\prime = AM$ with $M$ an invertible matrix. Replacing 
$A$ by $A^\prime$ and $c$ by $M^Tc$ yields the required reformulation, since
$
M^T c \, = \, M^T A^* \by  \, = \, A^{\prime *} \by,
$
so the first $k$ components of $M^Tc$ are zero. 
\qed

Note that we  slightly abuse terminology by calling  
(\ref{dreffeas}) a reformulation of \eref{d}: to obtain 
(\ref{dreffeas}) we not only reformulate \eref{d}, but also change the underlying cone. 

We now contrast Theorem \ref{sdlemma} with Theorem \ref{splemma}. 
In the former the minimal cone of $(P)$ can be written as 
$$
K \cap y_1^\perp \cap \dots \cap y_k^\perp,
$$
where $(y_1, \dots, y_k, y_{k+1})$ is some feasible solution in \eref{Dstrong}.
In the latter 
the minimal cone of the feasible set of 
\eref{d} is displayed by simply performing elementary row operations on 
the constraints. Thus, letting $G$ denote this minimal cone,
it is easy to convince a ``user'' that the Lagrange dual of 
\eref{dreffeas}, namely 
\beq
\ba{rrl}
\sup & \sum_{i=k+1}^m x_i c_i^\prime \\
& \sum_{i=1}^m x_i a_i^\prime & \! \leq_{G^*}  b 
\ena
\eeq
is an exact dual: the proof of this fact 
is some $y \in \ri G, \,$ which is  feasible in (\ref{dreffeas}), 
and the constraint set of (\ref{dreffeas}) 
which proves that all feasible solutions are in $G.$

\co{
To illustrate Theorem \ref{splemma}, we continue Example \ref{a1a2a3-ex}:
\bex{\rm (Example \ref{a1a2a3-ex} continued) The feasible set of this example, rewritten 
in an equality constrained form, looks like
$$
\ba{rcl} 
\bpx 0 & 0 & 0 \\ 0 & 0 & 0 \\ 0 & 0 & 1 & \epx \bullet y & = & 0 \\
\bpx 0 & 0 & -1 \\ 0 & 2 & 0 \\ -1 & 0 & 0 \epx \bullet y & = & 0 \\
\bpx 0 & 0 & 0 \\ 0 & 0 & 1 \\ 0 & 1 & 0 \epx \bullet y & = & 0 \\
\bpx 1 & 0 & 0 \\ 0 & 0 & 0 \\ 0 & 0 & 0 \epx \bullet y & = & 1,
\ena
$$
and the first two constraint matrices are in $\fr(\psd{3}).$ 
}
\eex
}

To illustrate Theorem \ref{splemma}, 
we continue Example \ref{a1a2a3-ex}:
\bex{\rm (Example \ref{a1a2a3-ex} continued) 
	We can rewrite the feasible set of this example
in an equality constrained form; note that 
if $z$ is a feasible slack in \eref{a1a2a3-sys}, then the 
equations 
\beq \label{rewrite} 
\ba{rcl} 
\bpx 0 & 0 & 0 \\ 0 & 0 & 0 \\ 0 & 0 & 1 & \epx \bullet z & = & 0 \\
\bpx 0 & 0 & -1 \\ 0 & 2 & 0 \\ -1 & 0 & 0 \epx \bullet z & = & 0
\ena
\eeq
must hold,  and  the constraint matrices in \eref{rewrite}  form a sequence in $\fr(\psd{3}).$ 
(Of course $z$ must satisfy two other linearly independent constraints as well,
which we do not show for brevity.)
}
\eex
\co{
for all $b$ objective functions \eref{dreffeas} has strong duality with its 
dual 
\beq
\ba{rcl}
\sup & \sum_{i=k+1}^m c_i^\prime x_i \\
s.t. & \sum_{i=1}^m a_i^\prime x_i & \leq_{(K^* \cap a_1^{\prime \perp} \cap \dots \cap a_k^{\prime \perp)^*}}
\ena
\eeq
}
\co{
\bth \label{d-strictfeas}
Suppose that \eref{d} is feasible. Then it has a strictly feasible reformulation 
\beq\label{dreffeas}\tag{$D_{\rm ref, feas}$}
\ba{rcl}
\la a_{i}^{\prime}, y\ra & = & 0\,(i=1,\ldots,k)\\
\la a_{i}^{\prime}, y\ra & = & c_{i}^{\prime}\,(i=k+1,\ldots,m)\\
y & \geq_G & 0 
\ena
\eeq
where 
$$
k \geq 0, \, (a_{1}^{\prime},\ldots,a_{k+1}^{\prime})\in\fr(K^*), \, G = (K^* \cap a_1^{\prime \perp} \cap \dots \cap a_k^{\prime \perp)^*}
$$ 
and for all $b$ objective functions \eref{dreffeas} has strong duality with its 
dual 
\beq
\ba{rcl}
\sup & \sum_{i=k+1}^m c_i^\prime x_i \\
s.t. & \sum_{i=1}^m a_i^\prime x_i & \leq_{(K^* \cap a_1^{\prime \perp} \cap \dots \cap a_k^{\prime \perp)^*}}
\ena
\eeq
\enth
}

To find the $y_i$ in Algorithm \ref{algo:FRA} one needs to solve a certain pair of 
{\em reducing conic linear programs}
(over $K$ and $K^*$). This task  may not be easier than solving the
problems \eref{p} and \eref{d}, 
since the primal reducing conic LP is strictly feasible, but its dual is not 
(see e.g. 
\cite[Lemma 1]{Pataki:13}).   
We know of two approaches to overcome this difficulty. The first approach by Cheung et al \cite{CheWolkSchurr:12}
is using a modified subproblem whose dual is also strictly feasible. 
The second, by Permenter and Parrilo in \cite{PerPar:14} is a 
``partial'' facial reduction algorithm for SDPs, where 
they solve linear programming approximations of the SDP subproblems. 

\section{Certificates  of infeasibility and weak infeasibility in conic linear programming}
\label{sec-catalog}

We now describe a collection of certificates of infeasibility and weak infeasibility of
\eref{p} and of \eref{d} below in Theorem \ref{mainfarkas}, which contains 
Theorem I.  
In Theorem \ref{mainfarkas} we state the results for the dual problem first, since we will mostly use these later on. 
We need to recall the definition of $\ell_K$ from \eref{defellK}. 
\bth \label{mainfarkas} 
When $K$ is a general closed, convex cone, the following hold: 
\benum 
\item \label{d-infeas} 
\eref{d} is infeasible, if and only if it has a reformulation 
\beq\label{dref2}\tag{$D_{\rm ref}$}
\ba{rcl}
\la a_{i}^{\prime}, y\ra & = & 0\,(i=1,\ldots,k)\\
\la a_{k+1}^{\prime}, y\ra & = & -1\\
\la a_{i}^{\prime}, y\ra & = & c_{i}^{\prime}\,(i=k+2,\ldots,m)\\
y & \geq_{K^*} & 0,
\ena
\eeq
where $(a_{1}^{\prime},\ldots,a_{k+1}^{\prime})\in\fr(K^*) $ and $0 \leq k \leq \ell_{K^*} - 1.$ 

\item \label{d-not-str} \eref{d} is not strongly infeasible, if and only if there is 
$(y_1, \dots, y_{\ell+1}) \in \fr(K)$, such that $0 \leq \ell \leq \ell_K -1 $ and 
\beqast
\ba{ccl}
A^*y_{i} &=& 0 \,(i=1,\ldots,\ell)\\
A^*y_{\ell+1} &=& c. 
\ena
\eeqast
\item \label{p-infeas} $(P)$ is infeasible, 
if and only if there is $(y_{1}, \dots, y_{k+1})\in\fr(K)$ such that $0 \leq k \leq \ell_K - 1$ and 
$$
\ba{rccccl}
A^*y_{i}&=&0,&b^*y_{i}&=&0\;(i=1,\ldots,k)\\
A^*y_{k+1}&=&0,&b^*y_{k+1}&=&-1.
\ena
$$
\item \label{p-not-str} $(P)$ is not strongly infeasible, if and only if it has a reformulation 
\beq\label{pref}\tag{$P_{\rm ref}$}
\sum_{i=1}^{m}x_{i}a_{i}^{\prime}\leq_{K} b^{\prime},
\eeq
where $(a_{1}^{\prime},\ldots,a_{\ell}^{\prime}, b^\prime)\in\fr(K^*)$ and 
$0 \leq \ell \leq \min \, \{ \, m, \, \ell_{K^*} - 1 \, \}.$ 
\eenum
The facial reduction sequences can be chosen to be pre-strict in all parts.
\enth
\qed

Parts \eref{d-infeas} through \eref{p-not-str} in Theorem \ref{mainfarkas} should be
read separately: the $k$ integers in parts \eref{d-infeas} and \eref{p-infeas}, 
the $\ell$ in parts \eref{d-not-str} and  \eref{p-not-str}, etc.  
may be different. 
We use the current notation for brevity. 
Also note that since $K$ is an arbitrary closed, convex cone, 
in the reformulations we only use elementary row operations 
(cf. Definition \ref{reform-def}). 

For the reader's sake we first discuss and illustrate Theorem \ref{mainfarkas} and prove a 
corollary, before 
proving Theorem \ref{mainfarkas} itself. 

We first look at the simplest cases, when $k=0$ or $\ell=0.$ 
We can choose $k=0$ in part (1) exactly if there is $x$ such that $Ax \geq_K 0$ and 
$\la c, x \ra = -1, \,$ i.e., if \eref{d} is strongly infeasible. 
Similarly, we can choose $\ell = 0$ in part (2) iff there is $y_1 \in K^*$ such that 
$A^* y_1 = c, \,$ i.e., if \eref{d} is actually feasible. 
\co{ 

We have
\beqast
{\rm We \, can \, choose \,} k=0 \, {\rm in \, part \eref{d-infeas}} & \LRA & \exists \, x \, s.t. \, Ax \in K, \, \la c, x \ra = - 1 \\
                                                                                                  & \LRA & \eref{d} \, {\rm is \, strongly \, infeasible}
\eeqast
Similarly, 
\beqast
{\rm We \, can \, choose \, }  \ell =0 \, {\rm in \, part \, } \eref{d-not-str}  & \LRA & \exists \, y_1 \in K^*  \, s.t. \, A^* y_1 = c  \\
& \LRA & \eref{d} \, {\rm is \, feasible}.
\eeqast}
Similar statements hold for parts \eref{p-infeas} and \eref{p-not-str}, and we leave the details to the reader. 

We next examine the ''easy'', i.e., the ''if'' directions. The ''if'' directions for 
parts \eref{d-infeas} and \eref{d-not-str} were proved after Theorem I. The ''if'' 
direction for part \eref{p-infeas} can be proved similarly to the proof of weak duality 
in Theorem \ref{sdlemma}. Next we prove the ''if'' direction in 
part (\ref{p-not-str}). It suffices to show that the subsystem 
\beq\label{george} 
\sum_{i=1}^{\ell}x_{i}a_{i}^{\prime}\leq_{K} b^{\prime},
\eeq
is not strongly infeasible, so suppose it is.
Applying the alternative system $(\palt)$  (defined in the Introduction)
to \eref{george}, we find that the system 
\beq \label{elaine} 
\ba{rcl} 
\la a_i^\prime, y \ra & = & 0 \, (i=1, \dots, \ell) \\
\la b^\prime, y \ra  & = & -1 \\
                       y & \geq_{K^*} & 0 
\ena
\eeq
is feasible. Since $(a_1^\prime, \dots, a_\ell^\prime, b^\prime)  \in \fr(K^*),$ we have
$$
b^\prime \in (K^* \cap a_1^{\prime \perp} \cap \dots \cap a_\ell^{\prime \perp})^*,
$$
thus for any feasible solution of (\ref{elaine}) we must have $\la b^\prime, y \ra \geq 0, \,$ 
a contradiction. 

Note the similarity between parts (1) and (4). 
Part (1) provides a certificate of infeasibility of \eref{d} in the form 
of a trivially infeasible subsystem. 
Part (4) gives an analogous certificate of ``not strong infeasibility'' of \eref{p} (since the subsystem 
\eref{george} has this status). 

All parts of Theorem \ref{mainfarkas} are useful. 
Part (1) allows us to generate all infeasible conic 
LP instances over cones, whose facial structure (and hence their facial reduction cone)
is well understood: to do so, we only need to generate systems of the form \eref{dref} and reformulate them.
For a more detailed discussion of the SDP case, we refer to Sections \ref{sec-catalog-sdp} 
and \ref{sec-generate}; 
for the case of conic LPs over {\em smooth cones}, which generalize $p$-order cones, we refer to Section 
\ref{sect-conclude}. 
By Part (4) we can systematically generate conic linear systems that are not strongly infeasible, though 
this seems less interesting. 

We will use parts (1) and (2) together to find our geometric corollaries
(in Section \ref{sec-geometry}), and to generate weakly infeasible SDPs (in Section 
\ref{sec-generate}). 

Part (3) strengthens the infeasibility certificate obtained by Waki and Muramatsu in 
\cite{WakiMura:12}: our certificate is essentially the same as theirs, the only difference is that 
we 	state it as a conic linear system. 
	
Part (4) is related to the recent paper of Lourenco et al \cite{Lourenco:13}. 
The authors there show that if a semidefinite system of the form \eref{p} 
is weakly infeasible, then a sequence 
$(a_1^\prime, \dots, a_\ell^\prime) \in \fr(\psd{n})$ 
can be found by taking linear combinations of 
the $a_i$ and applying rotations. (In fact, one can make 
$(a_1^\prime, \dots, a_\ell^\prime)$ to be a regularized facial reduction sequence, defined 
in Definition \ref{reg-def}.) 
In contrast, our results apply to general conic LPs; 
we  exactly characterize systems that are {\em infeasible} and systems that are 
{\em not strongly infeasible}. 
Putting these parts together 
yields our geometric corollaries 
(in Section \ref{sec-geometry}) and our algorithm to generate weakly infeasible 
SDPs (in Section \ref{sec-generate}). 

Example \ref{ex1} already illustrates parts \eref{d-infeas} and \eref{d-not-str}.
A larger example, which also depicts the frontier of $A^*K^*,$ follows;
in Example \ref{pic-ex} we choose $K = K^* = \psd{3}.$  

\begin{figure}[H]
\centering
\includegraphics[scale=0.5]{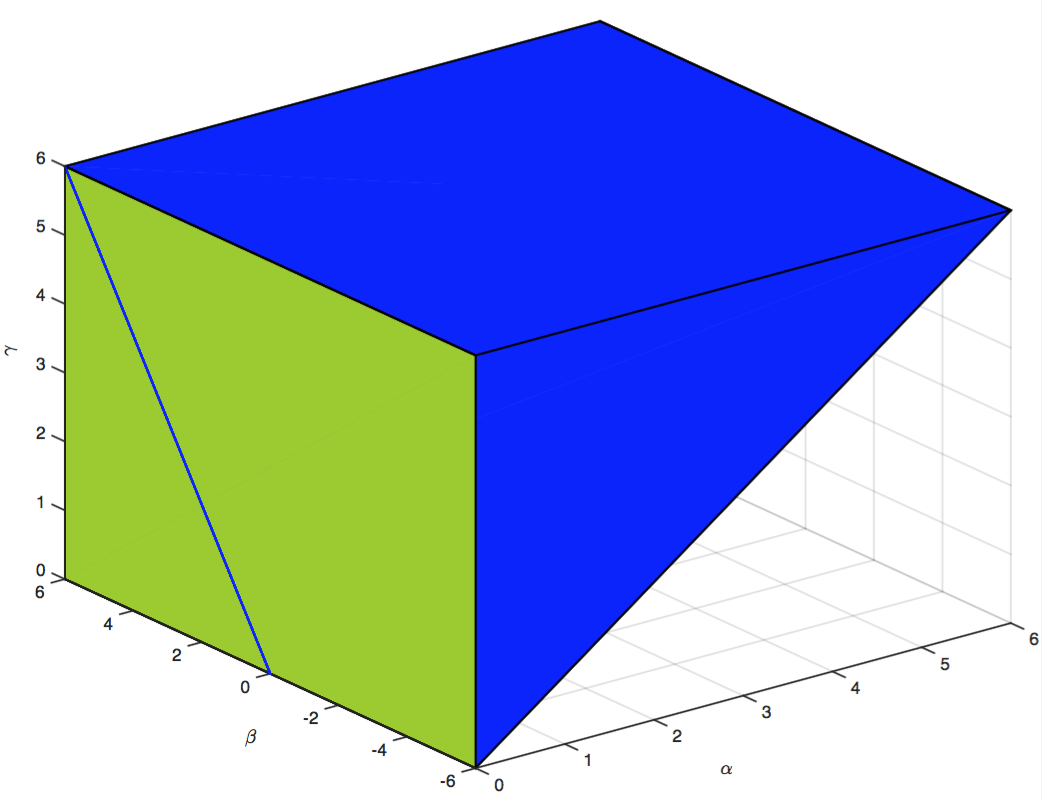} 
\caption{The set $A^*\psd{3}$ is in blue, and its frontier is in green} 
\label{frontier}
\end{figure}

\bex \label{pic-ex} {\rm Let 
\beq \label{ai-pic-ex} 
a_1 \, = \, \bpx 1 & 0 & 0 \\ 0 & 0 & 0 \\ 0 & 0 & 0 \epx, \, 
a_2 \, = \, \bpx 0 & 0 & 1 \\ 0 & 1 & 0 \\ 1 & 0 & 0 \epx, \, 
a_3 \, = \, \bpx 1 & 0 & 0 \\ 0 & 1 & 0 \\ 0 & 0 & 0 \epx.
\eeq 
Then we claim  
\beqa \label{cla} 
\cl (A^* \psd{3}) & = & \{ (\alpha, \beta, \gamma)^T \, : \, \gamma \geq \alpha \geq 0 \, \},  \\ \label{fronta} 
\front(A^* \psd{3}) \, = \, \cl (A^* \psd{3}) \setminus A^* \psd{3} & = & \{ (0, \beta, \gamma)^T \, | \, \beta \neq \gamma \geq 0 \, \}.
\eeqa
Indeed, the inclusion $\subseteq$ in \eref{cla} follows by calculating 
$A^* y$ for $y \succeq 0.$  
To see the reverse inclusion,  let $(\alpha, \beta, \gamma)^T$ be an element of the right hand side, and 
$$
y := \bpx
\alpha + \eps  & 0 & y_{13} \\
0        & \gamma - \alpha + \eps & 0 \\
y_{13} & 0                                  & y_{33}
\epx, \, 
$$
where $\eps >0, \, y_{13} \, = \, 0.5 (\beta - \gamma - \eps + \alpha)$, and $y_{33}$ is chosen to ensure $y \succeq 0. \,$ Then 
$$
A^* y = (\alpha+\eps, \beta, \gamma + 2 \eps)^T, 
$$
and letting $\eps \rightarrow 0$ completes the proof of $\supseteq.$ 
Equation \eref{fronta} follows by case checking.

The set $A^* \psd{3}$ is shown in Figure \ref{frontier} in blue, and its frontier in green.
Note that the blue diagonal segment inside the green frontier actually belongs to $A^*\psd{3}.$ (For better visibility, the $\beta$ axis goes from positive to negative 
 in Figure \ref{frontier}.)

To see how parts (1) and (2) of Theorem \ref{mainfarkas} certify that elements of 
$\front(A^* \psd{3})$ are indeed in this set, for concreteness, consider the system
\beq \label{noncl-ex} 
\ba{rcl}
A^*y  & = & (0,-2,1)^T \\
y      & \succeq & 0,
\ena
\eeq
which is weakly infeasible. 
The operations: 1)  multiply the second equation by $\frac{2}{3}$ and 2) add $\frac{1}{3}$ times the third equation to it, 
bring \eref{noncl-ex} into the form of \eref{dref} and show that it is infeasible. 
Using part (2) of Theorem \ref{mainfarkas}, the following $y_1$ and $y_2$  prove that it is not strongly infeasible:
\beq \label{y-ex-pic}
y_1 \, = \, \bpx 0 & \;0\; & 0 \\ 0 & \;0\; & 0 \\ 0 & \;0\; & 1 \epx, \, 
y_2 \, = \, \bpx 0 & \;0\; & -3/2 \\ 0 & \;1\; & 0 \\ -3/2 & \;0\; & 0 \epx,
\eeq
since $A^*y_1 = 0, \, A^*y_2 = (0,-2,1)^T, \, (y_1, y_2) \in \fr(\psd{3}).$ 
}
\eex
\co{
\compilepdf{
\begin{figure}[H]
\centering
\includegraphics[scale=0.6]{pic1}
\caption{The set $A^*\psd{3}$ is in blue, and its frontier is in green} 
\label{frontier}
\end{figure}
}}

To illustrate parts (3) and (4) of  Theorem \ref{mainfarkas}, 
we modify Example \ref{a1a2a3-ex} by simply exchanging two constraint matrices.
\bex \label{ex4} {\rm (Example \ref{a1a2a3-ex} continued) 
The  semidefinite system below is weakly infeasible.
\beq \label{a1a2a3-infeas}
\ba{rl}
   & x_1 \bpx 1 & 0 & 0 \\ 0 & 0 & 0 \\ 0 & 0 & 0 \epx + x_2 \bpx 0 & 1 & 0 \\ 1 & 0 & 0 \\ 0 & 0 & 0 \epx \preceq \bpx 0 & 0 & 1 \\ 0 & 1 & 0 \\ 1 & 0 & 0 \epx.
\ena
\eeq
To prove it is infeasible, we use part \eref{p-infeas} of Theorem \ref{mainfarkas} with 
$(y_1, y_2, y_3) \in \fr(\psd{3}), \, $ where $y_1, y_2$ are given in \eref{y1y2} and 
$$
y_3 = \bpx 0   & 0 & -1/2 \\
           0   & 0   & 0 \\
         -1/2  & 0   & 0 
       \epx. \, 
$$
To prove it is not strongly infeasible, we use part \eref{p-not-str}. We write 
$a_1, a_2,$ and $b$ for the constraint matrices, and observe that 
$(a_1, b) \in \fr(\psd{3}), $ and is pre-strict (we can also observe 
that $(a_1, a_2, b) \in \fr(\psd{3})). $ 
}
\eex

The next corollary states a ''coordinate-free'' version of Theorem \ref{mainfarkas} to address
a basic question in the theory of conic LPs: 
given a (weakly) infeasible 
conic linear system 
\beq \label{HK2} 
H \cap K,
\eeq
where $H$ is an affine subspace, what is the maximal/minimal dimension of an affine subspace $H^\prime$  with 
$H^\prime \supseteq H$ (or $H^\prime \subseteq H$) such that 
$H^\prime \cap K$ has the same feasibility status as \eref{HK2}? 
Note that by  (weak) infeasibility of \eref{HK2} we mean (weak) infeasibility of a representation 
in either the primal \eref{p}  or the dual \eref{d} form.

For instance, if $K$ is polyhedral, and \eref{HK2} is infeasible,
then by Farkas' lemma we can take $H^\prime$ as an affine subspace defined by a single equality constraint.
\co{
	contains $H$ (or is contained in $H$) so that replacing $H$ by $H^\prime$ does not change 
	the feasibility status of \eref{HK2}? } 

To further illustrate this question, 
\co{consider the semidefinite system 
\beq \label{kramer}
\ba{rcl} 
\bpx 1 & 0 & 0 \\ 0 & 0 & 0  \\ 0 & 0 & 0 \epx \bullet y & = & 0, \\
\bpx 0 & 0 & 1 \\ 0 & 1 & 0  \\ 1 & 0 & 0 \epx \bullet y & = & -1,
\ena
\eeq
which is weakly infeasible. Dropping the first constraint keeps it weakly infeasible,
and so does adding a constraint that fixes $y_{12}$ to zero. }
 let us revisit  Example \ref{ex4}. Here we can drop variable $x_2, \,$ or 
	add 
	$$
	x_3 \bpx 0 & 0 & 0 \\ 0 & 0 & 1 \\ 0 & 1 & 0 \epx
	$$ 
	to the left hand side, while keeping the system \eref{a1a2a3-infeas} weakly infeasible. 

Note that the {\em codimension} of a set $S$ is defined as the dimension of the underlying space minus the dimension of $S, \,$ and it is denoted by 
$\codim \, S.$ 

\bcor \label{dimension} 
The following hold.
\benum
\item If \eref{HK2} is infeasible, then there is $H^\prime \supseteq H$ such that 
$$
\codim H^\prime \leq \ell_K \, {\rm and} \, H^\prime \cap K \; {\rm is \; infeasible}.
$$
\item If \eref{HK2} is not strongly infeasible, then there is $H^{\prime \prime} \subseteq H$ such that 
$$
\dim H^{\prime \prime} \leq \ell_{K^*} -1 \, {\rm and} \, H^{\prime \prime} \cap K \; {\rm is \; not \; strongly \; infeasible}.
$$
\item If \eref{HK2} is weakly infeasible, then there is $H^{\prime \prime} \subseteq H \subseteq H^\prime$ as in parts (1) and (2) such that 
$$
H^{\prime} \cap K \; {\rm and} \; H^{\prime \prime} \cap K \; {\rm are \; both \;  weakly \; infeasible}.
$$
\co{\item If \eref{HK2} is weakly infeasible, then there is $H^{\prime \prime} \subseteq H \subseteq H^\prime$ as in parts (1) and (2) such that 
	$$
	\codim H^\prime \leq \ell_K \, \dim H^{\prime \prime} \leq \ell_K -1 \, \; {\rm and} \; H^{\prime} \cap K \; {\rm and} \; H^{\prime \prime} \cap K \; {\rm are \; both \;  weakly \; infeasible}.
	$$}

\eenum
\ecor 
\pf{} For part (1) we 
represent \eref{HK2} as a dual type problem \eref{d} (with $K$ in place of $K^*$), and 
apply part \eref{d-infeas} of Theorem \ref{mainfarkas}. We let
$H^\prime$ be the affine subspace defined by the first $k+1$ constraints in  \eref{dref}, and 
deduce
$$
\codim H^\prime \leq k+1 \leq \ell_K - 1 +1 = \ell_K,
$$
as required. 
For part (2) we represent \eref{HK2} as a primal type problem \eref{p} and apply 
part \eref{p-not-str} of Theorem \ref{mainfarkas}. We let 
$H^{\prime \prime} = \lin \{ a_1^\prime, \dots, a_\ell^\prime \} + b^\prime, \,$ where the 
$a_i^\prime$ and $b^\prime$ are given in \eref{pref}. 
Then \co{
Since $(a_1^\prime, \dots, a_\ell^\prime, b^\prime)$  is pre-strict, Theorem 
\ref{frlemma} implies that $a_1^\prime, \dots, a_\ell^\prime$ are linearly independent, so }
$$
\dim H^{\prime \prime} \leq \ell \leq \ell_{K^*} -1, 
$$
and this completes the proof. 

For part (3) we choose $H^\prime$ as in part (1). Since \eref{HK2} is 
not strongly infeasible,  and $H^\prime \supseteq H, \,$ the system $H^\prime \cap K$ is also not 
strongly infeasible; thus it is weakly infeasible.
We construct $H^{\prime \prime}$ as in part (2) with an analogous justification.
\qed

We now prove the ''only if'' parts in Theorem \ref{mainfarkas}.

\pf{of Theorem \ref{mainfarkas}:} 
We note in advance that the upper bounds on $k$ and $\ell$ will follow from the facial reduction 
sequences being pre-strict. 

To make the proofs concise, we prove the statements out of order, starting with \eref{p-infeas}.  
To see the ''only if'' part in part (3),  assume that $(P)$ is infeasible. and 
consider the conic LP 
\beq \label{supx0}
\sup\{x_{0}: Ax-bx_{0}\leq_{K}0\}, 
\eeq
which has value $0.$ 
The proof is complete by letting  
$(y_1, \dots, y_{\ell+1})$ to be a feasible solution to the exact dual of \eref{supx0}, 
given  in Theorem \ref{sdlemma}. 
\co{By pre-strictness, $y_1, \dots, y_\ell$ are linearly 
independent. Since $b^* y_{\ell+1} \neq 0, \,$  also 
$y_1, \dots, y_{\ell+1}$ are linearly independent, and this completes the proof. }

To prove  the ''only if'' part in \eref{d-infeas} we fix $\by$ such that 
$\{ \, y \, | \, A^* y = c \} = \by + \N(A^*).$ Choosing suitable generators for $\N(A^*)$ we can rewrite 
\eref{d} in the primal form, and apply part \eref{p-infeas}. Thus  there is 
$k \geq 0$ and a pre-strict $\, (a_{1}^{\prime},\ldots,a_{k}^{\prime}, a_{k+1}^\prime)\in\fr(K^*)$ such that
$$
\ba{rcll}
a_{i}^\prime      & \in   & \R(A) \cap \by^\perp &  ( i=1, \dots, k)  \\
a_{k+1}^\prime &  \in  & \R(A), &  \la a_{k+1}^\prime, \by \ra = -1. 
\ena
$$
\co{$$
\ba{rcll}
a_{i}^\prime   & \in & \R(A) \cap \by^\perp \,  (i=1, \dots, k), \\
a_{k+1}^\prime &  \in  & \R(A), \, \la a_{k+1}^\prime, \by \ra = -1. 
\ena
$$
} 
Since $\, (a_{1}^{\prime},\ldots,a_{k}^{\prime}, a_{k+1}^\prime)$ is pre-strict, 
$a_1^\prime, \dots, a_k^\prime$ are linearly independent. Since 
$\la a_{k+1}^\prime, \by \ra \neq \la a_{k}^\prime, \by \ra, $ also $a_1^\prime, \dots, a_k^\prime, a_{k+1}^\prime$ 
are linearly independent. The proof now can be completed verbatim as the proof of Theorem 
\ref{splemma}. 

For the ''only if'' part of 
statement \eref{d-not-str} we note that since \eref{d} is not strongly infeasible, the alternative 
system $(\dalt)$ is infeasible. We can view $(\dalt)$ as a conic linear system over the cone 
$K \ti \{0\}, \,$ i.e., 
\beq \label{Axc} 
\ba{rcl}
- Ax & \leq_K & 0 \\
\la -c,x \ra & \leq_{ \{0\}} & 1. 
\ena
\eeq
By part \eref{frkprop3} of Lemma \ref{frkprop}, and 
$\fr_k(\{0 \}) = \rad{k}, \,$  we find  that for all $k \geq 0$ 
\beq \nonumber 
\bigl( (y_1, z_1), \dots, (y_k, z_k) \bigr) \in \fr_k(K \ti \{0\}) \, \LRA \, (y_1, \dots, y_k) \in \fr_k(K) \, {\rm holds.} 
\eeq
Thus applying  part \eref{p-infeas} to the system \eref{Axc},  we deduce that there is 
$k \geq 0, \,$ and a  pre-strict 
$(y_{1},\ldots,y_{k+1}) \in \fr(K)$ and $(z_{1},\ldots,z_{k+1}) \in \rad{k+1}$ such that 
$$
\ba{rccccl}
A^*y_{i} + c z_i &=&0, &z_i & = & 0 \;(i=1,\ldots,k)\\
A^*y_{k+1} + c z_{k+1} &=&0, &z_{k+1} & = & -1,
\ena
$$
so our claim follows. 

Finally, to prove the ''only if'' part of \eref{p-not-str}, we note that since \eref{p} is not strongly infeasible, 
the system $(\palt)$ is infeasible, hence 
by part \eref{d-infeas} it has a reformulation 
\beq\label{paltref}
\ba{rcl}
\la a_{i}^{\prime}, y\ra &=& 0\,(i=1,\ldots,\ell)\\
\la b^{\prime}, y\ra &=& -1\\
\la a_{i}^{\prime}, y\ra &=& c_{i}^{\prime}\,(i=\ell+1,\ldots,m)\\
y & \geq_{K^*} & 0,  
\ena
\eeq
where $(a_1^\prime, \dots, a_\ell^\prime, b') \in \fr(K^*)$ for some $\ell \geq 0$ and 
the $c_i^\prime$ are suitable reals. (In (\ref{paltref}) we use $b^\prime$ to denote the 
left hand side of the special constraint, whose right hand side is $-1.$)  
Since in $(\palt)$ the only constraint with a nonzero right hand side is 
$\la b, y \ra = -1, \,$ we must have $b' = b + A \mu$ for some $\mu \in \rad{m}.$ 
Thus (\ref{pref}) is a reformulation of \eref{p} as required. 
\qed

\section{Certificates of infeasibility and weak infeasibility in  semidefinite programming} 
\label{sec-catalog-sdp}

In this section we specialize the certificates of infeasibility and weak infeasibility 
of Section \ref{sec-catalog} to semidefinite programming.
For this purpose we first introduce regularized facial reduction sequences in $\psd{n}:$ 
these sequences have a certain staircase like structure and 
we will use them in Theorem \ref{mainfarkas-sdp} 
in place of the usual facial reduction sequences in Theorem \ref{mainfarkas}. 

\begin{Definition} \label{reg-def} 
The set of order $k$ regularized facial reduction sequences for $\psd{n}$ is 
\beqast
\regfr_k(\psd{n}) & = & \biggl\{ (y_{1},\ldots,y_{k}): \, y_i \in \sym{n}, \, y_i \, = \, 
\bordermatrix{
& \overbrace{\qquad \qquad \qquad}^{\textstyle p_{1}+\ldots+p_{i-1}} & \overbrace{\qquad}^{\textstyle p_{i}} & \overbrace{\qquad \qquad \qquad\quad}^{\textstyle n-p_{1}-\ldots-p_{i}} \cr
& \times  &  \times  &  \times \cr
& \times  &  I   &  0 \cr
& \times  &  0  &  0 \cr} \\ 
          &  & \hspace{2.8cm} {\rm where} \, p_i \geq 0, \, i=1, \dots, k \biggr\}, \, 
\eeqast
where the $\times$ symbols correspond to blocks with arbitrary elements. 

We
drop the subscript  $k$ if its value is clear from the context or if it is irrelevant.
We say that $(y_1, \dots, y_k) \in \regfr_k(\psd{n})$ has block sizes $p_1, \dots, p_k$ 
if the order of the identity block in $y_i$ is $p_i$ for all $i.$ 
\end{Definition} 

Note that if $(y_1, \dots, y_k)  \in \regfr(\psd{n}) \,$ has block sizes $p_1, \dots, p_k, \,$ then it is strict, if and only if $p_i > 0$ for all $i.$ 

Also note that the constraint matrices in several examples, e.g. in Examples \ref{ex1} and  \ref{ex4},  actually form regularized facial reduction sequences. For example, the matrices in system 
(\ref{a1a2a3-infeas}) form a sequence in $\fr_3(\psd{3})$ with block sizes $1,0,1, \,$ respectively. 

Clearly,
$$
\regfr_k(\psd{n}) \, \subsetneq \fr_k(\psd{n}) \, 
$$
holds when  $n \geq 2.$  
However, Lemma \ref{tclemma} below shows that the set $\regfr_k(\psd{n})$ 
represents the larger set  $\fr_k(\psd{n}),$ since 
any element of $\fr_k(\psd{n})$ can be rotated to reside in $\regfr_k(\psd{n}).$ 
The proof of Lemma \ref{tclemma} is given in Appendix \ref{app-proofs}. 
\ble \label{tclemma}
Let $(y_1, \dots, y_{k}) \in \fr(\psd{n}).$ Then there is an invertible matrix $t$ such that 
$$
(t^Ty_1t, \dots, t^Ty_{k}t) \in \regfr(\psd{n}).
$$ 
\ele
\qed

The main result of this section follows.
\bth \label{mainfarkas-sdp} 
In the case of semidefinite programming, 
i.e., when $K = K^* = \psd{n},$ the following hold: 
\benum 
\item \label{d-infeas-sdp} 
\eref{d} is infeasible, if and only if it has a reformulation 
\beq\label{dref3}\tag{$D_{\rm ref, sdp}$}
\ba{rcl}
\la a_{i}^{\prime}, y\ra &=& 0\,(i=1,\ldots,k)\\
\la a_{k+1}^{\prime}, y\ra &=& -1\\
\la a_{i}^{\prime}, y\ra &=& c_{i}^{\prime}\,(i=k+2,\ldots,m)\\
y & \succeq & 0, 
\ena
\eeq
where $(a_{1}^{\prime},\ldots,a_{k+1}^{\prime})\in\regfr(\psd{n})$.
\item \label{d-not-str-sdp} \eref{d} is not strongly infeasible, if and only if
there is $(y_1, \dots, y_{\ell+1}) \in \regfr(\psd{n})$  such that
\beqast
\ba{ccl}
A^{\prime *}y_{i} &=& 0 \,(i=1,\ldots,\ell)\\
A^{\prime *}y_{\ell+1} &=& c^{\prime},
\ena
\eeqast
where $(A^{\prime}, c^\prime)$ is the data of a suitable reformulation of \eref{d}. 
\item \label{p-infeas-sdp} $(P)$ is infeasible, 
if and only if 
there is  $(y_{1}, \dots, y_{k+1}) \in \regfr(\psd{n})$ such that
$$
\ba{rccccl}
A^{\prime *} y_{i}&=&0,&b^{\prime *}y_{i}&=&0\;(i=1,\ldots,k)\\
A^{\prime *} y_{k+1}&=&0,&b^{\prime *}y_{k+1}&=&-1,
\ena
$$
where $(A^{\prime}, b^\prime)$ is the data of a suitable reformulation of \eref{p}.

\item \label{p-not-str-sdp} $(P)$ is not strongly infeasible, if and only if 
it has a reformulation 
\beq\label{prefsdp}\tag{$P_{\rm ref, sdp}$}
\sum_{i=1}^{m}x_{i}a_{i}^{\prime} \preceq  b^{\prime}
\eeq
where $(a_{1}^{\prime },\ldots,a_{\ell}^{\prime}, b^{\prime}) \in\regfr(\psd{n})$ for 
some $\ell \geq 0.$ 
\eenum 
In all parts the facial reduction sequences can be chosen as pre-strict, and $k$ and $\ell$ at
most $n-1.$ 

\co{Furthermore,  $a_1^\prime, \dots, a_{k+1}^\prime$ can be chosen to be 
linearly independent in part (1); and $y_{1}, \dots, y_{\ell+1}$ can be chosen to be 
linearly independent in part (4).}
\qed
\enth

Again, all parts of Theorem \ref{mainfarkas-sdp} should be read 
separately, i.e., the reformulations, and the $k$ and $\ell$ integers in parts (1) through (4) 
may all be different. 

\co{Note that Theorem \ref{mainfarkas-sdp} differs from Theorem \ref{mainfarkas} in three ways.
First, Theorem \ref{mainfarkas-sdp} uses 
regularized facial reduction sequences (see Definition \ref{reg-def}) in place of the usual facial reduction sequences (see Definition \ref{frk-def}). 
Second, when ''reformulating'' \eref{p} or \eref{d} in Theorem \ref{mainfarkas-sdp} 
we are allowed to choose a suitable $T \in \aut(\psd{n}) \,$ 
(see (\ref{def-autk})) 
and to replace all $a_i$ by $T a_i$ and  $b$ by $Tb \,$ --- this is because $\psd{n}$ is self-dual.
(It is well-known that $T \in \aut(\psd{n})$ if and only if there is 
an invertible matrix $t$ such that 
$T(x) = t^T x t$ for $x \in \sym{n}.$) 

Third, in parts (2) and (3) of Theorem 
\ref{mainfarkas-sdp} we refer to a suitable reformulation of \eref{d} and of \eref{p}, while 
in parts (2) and (3) of Theorem \ref{mainfarkas} we simply refer to \eref{d} and \eref{p}. 
}

We emphasize three key differences between 
Theorem \ref{mainfarkas-sdp} and Theorem \ref{mainfarkas}.
	First, Theorem \ref{mainfarkas-sdp} uses 
	regularized facial reduction sequences  in place of the usual facial reduction sequences. 
	Second, when ''reformulating'' \eref{p} or \eref{d} in Theorem \ref{mainfarkas-sdp} 
	we are allowed to choose a suitable invertible matrix $t$ and replace 
	all $a_i$ by $t^T a_i t$ and $b$ by $t^Tbt$. This is because $\psd{n}$ is self-dual, and 
	the description of the automorphism group of $\psd{n}$ after equation (\ref{def-autk}).
	Third, in parts (2) and (3) of Theorem 
	\ref{mainfarkas-sdp} we refer to a suitable reformulation of \eref{d} and of \eref{p}, while 
	in parts (2) and (3) of Theorem \ref{mainfarkas} we simply refer to \eref{d} and \eref{p}. 
	
	Also, simply applying Theorem \ref{mainfarkas} with $K = K^* = \psd{n}$ we would get the upper bound 
	$\ell_{\psd{n}}-1 = (n+1) - 1 = n$ on $k$ and $\ell, \,$ 
		while Theorem \ref{mainfarkas-sdp} has the slightly stronger bound 
		$n-1.$ 

We note that part (1) in Theorem \ref{mainfarkas-sdp} recovers Theorem 1 in 
\cite{LiuPataki:15}. 

What is the main advantage of Theorem \ref{mainfarkas-sdp} vs. Theorem 
\ref{mainfarkas}? In Theorem \ref{mainfarkas-sdp} the ''easy,'' i.e., the ''if'' directions can be proved 
relying only on elementary 
linear algebra. For instance, to prove the ''if'' direction 
in part (1), let us assume that $y \succeq 0$ is 
feasible in \eref{dref3}, and that 
$a_1^\prime, \dots, a_k^\prime$ have block sizes $r_1, \dots, r_k, \,$ respectively.
Then $a_1^\prime \bullet y = 0, \,$ and $y \succeq 0$ imply that 
the first $r_1$ rows and columns of $y$ are zero. We then  inductively prove 
that the first $r_1 + \dots + r_k$ rows and columns of $y$ are zero, 
hence $a_{k+1}^\prime \bullet y \geq 0, \,$ which is a contradiction. For more  details, we refer to \cite{LiuPataki:15}. 
We can similarly prove  the ''if'' 
directions in parts (2) through (4).

Again, for the reader's sake we first illustrate Theorem \ref{mainfarkas-sdp}, and 
specialize Corollary \ref{dimension}.

\bex{\rm (Example \ref{ex1} continued)
	We revisit  the weakly infeasible instance 
	\beq \label{ex1-sys-2}
	\ba{rcl} 
	\bpx 1 & 0 \\ 0 & 0 \epx \bullet y & = & 0 \\
	\bpx 0 & 1 \\ 1 & 0  \epx \bullet y & = & -1 \\
	y & \succeq & 0.
	\ena
	\eeq
	Writing $a_1$ and $a_2$ for the constraint matrices, we see  that 
	$(a_1, a_2) \in \regfr(\psd{2}), \,$ so  this example illustrates part (1) of Theorem 
	\ref{mainfarkas-sdp}. 
Next, let 
$$
t = \bpx 0 & 1 \\ 1 & 0 \epx,
$$
and apply the transformation 
$t^T()t$ to $a_1$ and $a_2$ to obtain the system 
\beq \label{ex1-sys-3}
\ba{rcl} 
\bpx 0 & 0 \\ 0 & 1 \epx \bullet y & = & 0 \\
\bpx 0 & 1 \\ 1 & 0  \epx \bullet y & = & -1 \\
y & \succeq & 0. 
\ena
\eeq
This latter system illustrates part (2) of Theorem \ref{mainfarkas-sdp}, since
the fact that it is not strongly infeasible is proved by 
$$
y_1 = \bpx 1 & 0 \\ 0 & 0 \epx, \, y_2 = \bpx 0 & -1/2 \\ -1/2 & 0 \epx,
$$
and $(y_1, y_2) \in \regfr(\psd{2}).$ 
}
\eex

\bex{\rm (Example \ref{pic-ex} continued) The system \eref{noncl-ex}  in Example \ref{pic-ex} illustrates part (2) in Theorem \ref{mainfarkas-sdp}, as follows: let us exchange the first row with the third row and the first column with the third column in all 
	constraint matrices in (\ref{noncl-ex}) 
	to obtain a system which is also weakly infeasible.
	
	Let us denote this system by 
	\beq \label{noncl-ex2} 
	\ba{rcl}
	A^{\prime *} y  & = & (0,-2,1)^T \\
	y      & \succeq & 0.
	\ena
	\eeq
	Letting 
	\beq \label{y-ex-pic-2}
	y_1 \, = \, \bpx 1 & \;0\; & 0 \\ 0 & \;0\; & 0 \\ 0 & \;0\; & 0 \epx, \, 
	y_2 \, = \, \bpx 0 & 0 & -3/2 \\ 0 & 1 & 0 \\ -3/2 & 0 & 0 \epx,
	\eeq
	we have $A^{\prime *} y_1 = 0, \, A^{\prime *} y_2 = (0, -2, 1)^T, \, (y_1, y_2) \in \regfr(\psd{3}),$
	thus the system (\ref{noncl-ex2}) and $(y_1, y_2)$ are as required in part (2) of Theorem 
	\ref{mainfarkas-sdp}. 
}
\eex

\co{(The reader may n If we apply the same rotation to the 
$y_j$ matrices in equation \eref{y-ex-pic} (and for brevity, we write $y_j$ for $t^T y_j t$) 
we obtain )}
The reader can also check that 
Example \ref{ex4} illustrates parts (3) and (4) of Theorem \ref{mainfarkas-sdp}. 

\bcor
Consider the conic linear system \eref{HK2} and assume that the 
underlying cone is $K = \psd{n}.$ 
All statements of Corollary \ref{dimension} applied to this 
system hold with $\ell_{\psd{n}} - 1 = n$ in place of 
$\ell_{\psd{n}} = n+1.$ 

In other words we can choose $H^\prime$ and $H^{\prime \prime}$ such that 
$\codim H^\prime \leq n, \,$ and $\dim H^{\prime \prime}  \leq n-1.$ 
\ecor
\pf{} The proof is identical to the proof of 
Corollary \ref{dimension}, except we need to invoke 
Theorem \ref{mainfarkas-sdp} in place of Theorem 
\ref{mainfarkas}.  
\qed

\pf{of the ''only if'' parts in Theorem \ref{mainfarkas-sdp}} 
To see part \eref{d-infeas-sdp} we consider the reformulation given in 
part  \eref{d-infeas} of Theorem \ref{mainfarkas}, and a $t$ invertible matrix 
such that $$(t^T a_1^\prime t, \dots, t^T a_{k+1}^\prime t) \in \regfr(\psd{n}).$$ 
We replace $a_i^\prime$ by $t^T a_i^\prime t$ for all $i$ and obtain \eref{dref3}. 

To prove $k \leq n -1 \,$ in part (1), we note that the bound $k \leq n$ follows from 
part (1) of Theorem \ref{mainfarkas}, since $\ell_{\psd{n}} = n+1.$ 
If $k < n, \,$ then there is nothing to prove, so suppose 
$k = n. \,$ Since $(a_1^\prime, \dots, a_n^\prime)$ is strict, 
the block sizes of all the $a_i^\prime$ must be all $1,$ so the lower 
right  
$2 \ti 2$ blocks of $a_{n-1}^\prime$ and $a_n^\prime$ look like 
$$
\bpx 1 & 0 \\ 0 & 0 \epx  \, {\rm and} \, \bpx \ti & \ti \\ \ti & 1 \epx,
$$
where the $\ti$ symbols stand for arbitrary components.

Let us choose a suitable $\lambda$ and set  
$a_{n-1}^{\prime \prime} := a_{n-1}^\prime + \lambda a_n^\prime,$ so 
the lower right $2 \ti 2$ block of $a_{n-1}^{\prime \prime}$ is positive 
definite. It is then easy to check that 
$(a_1^\prime, \dots, a_{n-2}^\prime, a_{n-1}^{\prime \prime}, a_{n+1}^\prime) \in \fr(\psd{n}),$ and it is pre-strict. 
Finally, using Lemma \ref{tclemma},  
we pick a suitable invertible matrix $t$ and apply the rotation 
$t^T()t$ to this sequence, so the result is in $\regfr(\psd{n}). \,$ This completes the proof.

To see \eref{d-not-str-sdp} we consider the sequence 
$(y_1, \dots, y_{\ell+1}) \in \fr(\psd{n})$ 
given by part  \eref{d-not-str} of Theorem \ref{mainfarkas}, and a $t$ invertible matrix 
such that 
$$
(t^Ty_1t, \dots, t^Ty_{\ell+1}t) \in \regfr(\psd{n}).
$$ 
For $i=1, \dots, m$ and $j = 1, \dots, \ell+1$ we have 
\beq \nonumber 
\ba{rcl}
\la t^{-1}a_it^{-T}, t^T y_j t \ra & = & \la a_i, y_j \ra.
\ena
\eeq
We set $a_i^{\prime} := t^{-1}a_it^{-T}$ and replace $y_j$ by $t^T y_j t$ for all $i$ and $j$. 
The upper bound on $\ell$ follows by a similar argument as the upper bound on $k$ 
in part (1), and this completes the proof. 

The proof of \eref{p-infeas-sdp} is analogous to the proof of \eref{d-not-str-sdp}; 
and the proof of \eref{p-not-str-sdp} to the proof of \eref{d-infeas-sdp}, hence we omit these.
\qed

\section{Geometric corollaries} 
\label{sec-geometry}

In this section we  use the preceding results to address several 
basic questions in convex analysis. We begin by asking the  question:
\begin{center} 
Under what conditions is the linear image of a closed convex cone closed?
\end{center} 
This question is fundamental, due to its role in constraint 
qualifications in convex programming. Due to its importance, Chapter 9 in Rockafellar's classic text 
\cite{Rockafellar:70} is entirely 
devoted to it: see e.g. Theorem 9.1 therein. 
Surprisingly, the literature on the  subject (beyond \cite{Rockafellar:70} and other textbooks)
appears to be scant. Bauschke and Borwein \cite{BausBor:99}
gave a necessary and sufficient condition for the continuous image of a closed convex 
cone to be closed. Their condition  (due to its greater generality) is more involved 
than Theorem 9.1 in \cite{Rockafellar:70}. See also \cite{Aus:96}, and the references in 
\cite{Pataki:07}. We also refer to Bertsekas and Tseng 
for a study of a more general problem, when the intersection of a nested sequence of sets is nonempty. See Borwein and Moors \cite{BorweinMoors:09, BorweinMoors:10} for 
proofs that the set of linear maps under which the image is {not closed} is small both 
in terms of measure and category. 

For convenience we restate our question in an equivalent form: 
\begin{center} 
Given $A$ and $K, \,$ when is $A^*K^*$ closed?
\end{center} 
In \cite{Pataki:07} we gave the very simple {\em necessary} condition 
\beq \label{radir}
\R(A) \cap (\cl \dir(z, K) \setminus \dir(z,K)) = \emp
\eeq
for $A^*K^*$ to be closed: here 
$z$ is in the relative interior of $\R(A) \cap K, \,$ and 
$$
\dir(z,K) = \{ \, y \, | \, z + \eps y \in K \, {\rm for \, some \, } \eps > 0 \, \}$$ 
is the set of 
feasible directions at $z$ in $K.$ Note that \eref{radir} 
subsumes two seemingly unrelated classical
sufficient conditions for the closedness of $A^*K^*,$ as it trivially holds when $K$ is polyhedral, or when $z \in \ri K.$ 
It is also sufficient, when the set  
$$
K^* + F^\perp
$$
is closed, where $F$ is the minimal cone of $\R(A) \cap K.$ Thus \eref{radir} becomes 
an exact characterization when $K^* + F^\perp$ is closed for {\em all} $F$ faces of 
$K.$ Such cones are called {\em nice}, and reassuringly, 
most cones that occur in optimization (such as polyhedral, semidefinite, and $p$-order cones) are nice. 
Nice cones also play a role in simplifying constraint qualifications in conic LPs: see 
\cite{BorWolk:81,BorWolk:81B}. 

As a byproduct of the preceding results, here we 
obtain an {\em exact} and simple characterization of when $A^*K^*$ is closed,
 when $K$ is an arbitrary 
closed convex cone. 

We build on the following basic fact:
\beq
A^*K^* \; {\rm is \; not \; closed \;} \LRA \; \eref{d} \; {\rm is \; weakly \; infeasible \; for \; some \; } c.
\eeq
\bth \label{not-closed} The set $A^*K^*$ is not closed, if and only if there is 
$(a_1, \dots, a_{k+1}) \in  \fr_{k+1}(K^*)$ with $k \geq 1,$ and 
$(y_1, \dots, y_{\ell+1}) \in  \fr_{\ell+1}(K)$ with $\ell \geq 1$ such that 
\beq \label{aiyj-1}
\ba{rcl}
a_i & \in & \R(A) \, (i=1, \dots, k+1), \\
y_j & \in & \N(A^*) \, (j = 1, \dots, \ell)
\ena
\eeq
and 
\beq \label{aiyj-2} 
\ba{rcl}
\la a_i,  y_{\ell+1} \ra \, = \, \left\{ \ba{rl} 0 & {\rm if} \, i \leq k \\
                                                -1 & {\rm if} \, i = k+1. 
                                \ena \right.
\ena
\eeq
\enth
\pf{} 
Starting  with the forward implication, we choose $c$ such that \eref{d} is weakly infeasible.
We choose $a_1, \dots, a_k, a_{k+1}$ as the left hand side vectors in 
the first $k+1$ constraints in a reformulation of the form \eref{dref}, 
	which proves that \eref{d} is 
	infeasible; 
	and $(y_1, \dots, y_{\ell+1})$ to prove that it is not strongly infeasible: 
cf. parts \eref{d-infeas} and \eref{d-not-str}  in Theorem \ref{mainfarkas}.

For the backward implication, fix $a := (a_1, \dots, a_{k+1})$ and $y := (y_1, \dots, y_{\ell+1})$ 
as stated. First we prove that they can be assumed to be pre-strict, so 
suppose that, say, $a$ is not. Then for some $i < k$ we have 
$
a_{i+1} \in (K \cap a_1^\perp \cap \dots \cap a_i^\perp)^\perp, 
$
hence 
$$
K \cap a_1^\perp \cap \dots \cap a_i^\perp \, = \, K \cap a_1^\perp \cap \dots \cap a_i^\perp \cap a_{i+1}^\perp,
$$
so we can drop $a_{i+1}$ from $a$ while keeping all required properties of $a$ and $y.$ 
Continuing like this we arrive at both $a$ and $y$ being pre-strict, and to ease notation, we 
still assume $a \in \fr_{k+1}(K^*)$ and $y \in \fr_{\ell+1}(K).$ 

Now $a_1, \dots, a_k$ are linearly independent (by Theorem \ref{frlemma}). Since 
$\la a_{k+1}, y_{\ell+1} \ra \neq \la a_{k}, y_{\ell+1} \ra, $ so are $a_1, \dots, a_k, a_{k+1}.$ 

Thus we can expand $a_1, \dots, a_k, a_{k+1}$ into 
$$
A^\prime = [a_1, \dots, a_k, a_{k+1}, a_{k+2}, \dots, a_{m}],  \,\,  {\rm a \,\, basis \,\, of \,\,} \R(A), 
$$
and let 
$$
c^\prime = (0, \dots, 0, -1, \la a_{k+1}, y_{\ell+1} \ra, \dots, \la a_{m}, y_{\ell+1} \ra)^T.
$$
Write $A^\prime = AM, \,$ with $M$ an $m \ti m$ invertible matrix, and let
 $c = M^{-T}c^\prime.$ 
Then \eref{d} with this $c$ is weakly infeasible (since 
$A^{\prime *}= M^TA^*, \, c^\prime = M^Tc, \,$ so the reformulation with data 
$(A^\prime, c^\prime)$ proves that \eref{d} is infeasible;
and $(y_1, \dots, y_{\ell+1})$ proves  that  it is not 
strongly infeasible: cf. Theorem \ref{mainfarkas}). 
\qed

It is also of interest to characterize nice cones. To review previous results on nice cones we
recall that $y \in K^*$ is said to 
{\em expose} the set $G$ if $G \subseteq K \cap y^\perp; $ 
 a face $G$ of $K$ is said to be {\em exposed}, 
if it equals $K \cap y^\perp$ for some $y \in K^*;$ and it is not exposed iff 
\beq \label{FGperp}
K^* \cap G^\perp = K^* \cap F^\perp
\eeq
for some $F$ face of $K$ that strictly contains $G$ 
(note that the set of vectors that expose $G \subseteq K$ is $K^* \cap G^\perp, \,$ 
so (\ref{FGperp}) means that all vectors that expose $G$ actually expose the larger face $F$). Note that if $y \in K^*, \,$ then $K \cap y^\perp$ is a supporting 
hyperplane of $K.$
The cone $K$ is said to be {\em facially exposed, } if all of its faces are exposed.

Example \ref{nonexp} below shows a cone which is not facially exposed.
\bex \label{nonexp} {\rm Define $K$ as the sum of $\psd{2}$ and the cone comprising all nonnegative multiples of the
matrix 
$$
\bpx 0 & 1 \\ 
     1 & 0 \epx.
$$
\compilepdf{
\begin{figure}[H]
\begin{center}
\includegraphics[width=4in]{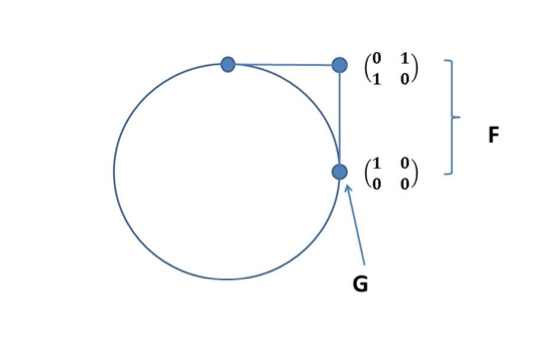}   
\caption{\small Cross section of a cone which is not facially exposed \normalsize}
\label{fig-nonexp}
\end{center}
\end{figure}
}
The cross-section of this cone is shown on Figure \ref{fig-nonexp}, with faces 
$G \subsetneq F$ that satisfy \eref{FGperp}.
 The face $G$ is not exposed, while $F$ is. 
}
\eex
For characterizations of nice cones, and a proof that they must be facially exposed, 
we refer to Pataki \cite{Pataki:12}; for an example of a facially exposed, but not nice cone, 
see Roshchina \cite{Vera:13}; and Chua and Tun\c{c}el
\cite{ChuaTuncel:08} for a proof that the linear pre-image of a nice cone is also nice.

Theorem \ref{not-closed} also leads to a characterization of when a cone is (not) nice: 
\bth \label{not-nice} 
Let $F$ be a face of $K.$ Then $K^* + F^\perp$ is not closed, if and only if 
there is 
$(a_1, \dots, a_{k+1}) \in  \fr_{k+1}(K^*)$ with $k \geq 1,$ and 
$(y_1, \dots, y_{\ell+1}) \in  \fr_{\ell+1}(K)$ with $\ell \geq 1$ such that 
\beq \label{aiyj-nice-1} 
\ba{rcl}
a_i & \in & \lin F  \, (i=1, \dots, k+1), \\
y_j & \in & F^\perp \, (j = 1, \dots, \ell)
\ena
\eeq
and 
\beq \label{aiyj-nice-2} 
\ba{rcl}
\la a_i,  y_{\ell+1} \ra \, = \, \left\{ \ba{rl} 0 & {\rm if} \, i \leq k \\
                                                -1 & {\rm if} \, i = k+1. 
                                \ena \right.
\ena
\eeq
\enth
\pf{} The result follows from Theorem \ref{not-closed} by considering a 
linear operator $A$ with $\R(A) = \lin F, \, \N(A^*) = F^\perp$ and 
noting that $A^*K^*$ is  not closed, iff $\N(A^*) + K^*$ is not closed 
(see Lemma 3.1 in \cite{Berman:73}). 
\qed

Theorems \ref{not-closed} and \ref{not-nice} provide a hierarchy of conditions, and 
it is natural to ask, how these relate to the simpler, but less general known conditions on closedness, and niceness. 
To address this  question, we need two definitions: 
\begin{Definition}
We say that the nonclosedness of $A^*K^*$ (of $K^* + F^\perp$) has a $(k+1, \ell+1)$-proof, if 
there is $(a_1, \dots, a_{k+1})$ and $(y_1, \dots, y_{\ell+1})$ as in Theorem 
\ref{not-closed} (Theorem \ref{not-nice}).
\end{Definition}

\begin{Definition} \label{def-sing}
The singularity degree of the system $H \cap K, \,$ where $H$ is an affine subspace,
is the minimum number of facial reduction steps needed by 
Algorithm \ref{algo:FRA} to find its minimal cone. 
\end{Definition}

For example, the singularity degree of 
the system \eref{a1a2a3-sys} is two: the facial reduction sequence given in 
\eref{y1y2} is shortest possible. 

\bth \label{nonclosed} The following hold:
\benum
\item Suppose that condition \eref{radir} is violated, and let $\ell$ be the 
singularity degree of $\R(A) \cap K.$ Then there is a $(2, \ell+1)$-proof of the nonclosedness of $A^*K^*.$ 
\item Suppose that $K$ has a nonexposed face, say $G, \,$ 
and $F$ is the smallest exposed face of $K$ that contains it. Then there is a $(2,2)$-proof that 
$K^* + F^\perp$ is not closed.
\eenum
\enth
Since the proof of this result is somewhat technical, we defer it to Appendix \ref{app-nonclosed}. 
It is also natural to ask, 
as to what values of $k$ and $\ell$ are actually necessary to prove nonclosedness of $A^*K^*$ (or of 
$K^* + F^\perp$). 
A recent result of Drusviyatsky et al \cite{Dima:15} shows a surprising connection 
between the degree of singularity of the dual problem \eref{d}, and 
the exposedness of the smallest face of $A^*K^*$ that contains $c$. 
For an equivalent, independently obtained  result for semidefinite programs,
see Gortler and Thurston \cite{GortlerThurston:14}. 

It would also be interesting to 
explore the connection to this result and we will do so in a followup paper. 

Also, in recent work, Roshchina and Tun\c{c}el gave a condition to 
strengthen the facial exposedness condition of \cite{Pataki:12}: it would be interesting to see
how their condition fits into our hierarchy. 

\section{Generating infeasible, and weakly infeasible SDPs}
\label{sec-generate}

We now turn to a practical aspect of our work, generating infeasible, and weakly infeasible SDP instances in the dual form 
\eref{d}. 
Having a library of such instances is important, 
since detecting infeasibility is a weak point of 
commercial and research codes: when they report that 
\eref{d} is infeasible,
they also return a certificate of infeasibility, namely a feasible solution to 
the alternative system $(\dalt).$ If \eref{d} is weakly infeasible, then 
$(\dalt)$ is also infeasible, so the 
returned certificate is necessarily inaccurate.

We first state an elementary algorithm, based on part (1) 
of Theorem \ref{mainfarkas-sdp}, to generate infeasible SDPs in the standard 
form (\ref{dref3}). 
\begin{algorithm}[H] 
\caption{Infeasible SDP}
\label{algo:non_infeas_SDP}
\begin{algorithmic}
\State 1: Choose integers $m, n, k, p_1, \dots, p_k > 0$ and $p_{k+1} \geq 0$ s.t. $k+1 \leq m, \, \sum_{i=1}^{k+1}p_{i}\leq n.$
\State 2: Let $(a_1, \dots, a_{k+1}) \in \regfr_{k+1}(\psd{n})$ with block sizes $p_1, \dots, p_{k+1}$ and $c_1 = \dots = c_k = 0, \, c_{k+1} = -1.$
\State 3: Let $a_{k+2}, \dots, a_m \in \sym{n}$ and $c_{k+2}, \dots, c_m \in \rad{}$ be arbitrary.
\end{algorithmic}
\end{algorithm}

By part (1) of Theorem \ref{mainfarkas-sdp} all 
infeasible SDPs (in the dual form \eref{d})  are obtained by taking  
some output of Algorithm \ref{algo:non_infeas_SDP}, and arbitrarily reformulating it. 
Algorithm \ref{algo:non_infeas_SDP} may generate a  strongly or a weakly infeasible SDP, and the latter outcome is 
likelier if $m$ is small with respect to $k+1 \,$
(in this case $(\dalt)$ is tightly constrained, so it is likely to be infeasible).

However, weak infeasibility of the output of 
Algorithm  \ref{algo:non_infeas_SDP}  is not guaranteed. 

Next we turn to generating weakly infeasible SDP instances with a proof  of weak 
infeasibility. 
In contrast to Waki's instances 
in \cite{Waki:12}, 
we will generate our instances by solving
simple systems of equations. In fact, we will define a natural class of 
weakly infeasible SDPs, and show that a simple algorithm generates {\em all}
instances in this class. The instances will be in a form so that one can 
verify  their weak infeasibility by an elementary linear algebraic argument.

Although our framework is different -- since we generate objects in an uncountably infinite set --
our algorithms to generate {\em all} SDP instances in a certain class 
fit into the framework of {\em listing} combinatorial objects, as cycles, paths, spanning trees and cuts:
see e.g., \cite{ReadTarj:75, ProvanShier:96}. 

We will use part (1) of Theorem \ref{mainfarkas} to find an infeasible instance,
and part (2) to find a $(y_j)$ sequence to prove that it is not strongly infeasible, so 
we will solve a {\em bilinear} system of equations over the 
$a_i$ and $y_j.$
While this may be difficult in general, it is easy if we impose a structure:  
we will require that the $(a_i)$ sequence be regularized (cf. Definition 
\ref{reg-def}), and that the $(y_j)$ have the same structure, but ``reversed'' in the sense defined 
below:

\begin{Definition} \label{revreg-def} 
The set of order $\ell$ reversed regularized facial reduction sequences in $\psd{n}$ is 
\beqast
\revregfr_\ell(\psd{n}) & = & \biggl\{ (y_{1},\ldots,y_{\ell}): \, y_j \in \sym{n}, \, y_j \, = \,
\bordermatrix{
& \overbrace{\qquad \qquad \qquad}^{\textstyle n-q_{1}-\ldots-q_{j}} & \overbrace{\qquad}^{\textstyle q_{j}} & \overbrace{\qquad \qquad \qquad\quad}^{\textstyle q_{1}+\ldots+q_{j-1}} \cr
& 0       &  0   &  \times \cr
& 0       &  I   &  \times \cr
& \times  &  \times  &  \times \cr} \\ 
          &  & \hspace{2.8cm} {\rm where} \, q_j \geq 0, \, j=1, \dots, \ell \biggr\}, \, 
\eeqast
where the $\times$ symbols correspond to blocks with arbitrary elements. 

We drop the subscript $\ell$, if its value is clear from the context or if it is 
irrelevant. 
We also say that $(y_{1},\ldots,y_{\ell}) \in \revregfr_\ell(\psd{n})$ has block sizes 
$q_1, \dots, q_\ell$ if the  order of the identity block in $y_j$ is $q_j$ for all $j.$ 
\end{Definition}
For instance, the $y_j$ matrices in Example \ref{ex1} are in 
$\revregfr(\psd{2}).$ 
\begin{Definition}
An SDP instance 
\beq \label{nonover} 
\ba{rcl}
A^*y & = & c \\
   y & \succeq & 0 
\ena
\eeq
is {\em nonoverlapping weakly infeasible}, if 
\benum
\item it is in the form \eref{dref} as in part (1) of Theorem \ref{mainfarkas} with 
$(a_1, \dots, a_{k+1}) \in \regfr(\psd{n})$ and 
$c_1 = \dots = c_k = 0, \, c_{k+1} = -1 \,$ for some $k \geq 1.$ 
\item There is $(y_1, \dots, y_{\ell+1}) \in \revregfr(\psd{n})$ 
as in part (2) of Theorem \ref{mainfarkas} which proves that it is not strongly infeasible; 
\item The block sizes $p_i$ of $(a_1, \dots, a_{k+1})$ and the block sizes $q_j$ of $(y_1, \dots, y_{\ell+1})$ satisfy 
\beq \label{piqj} 
\sum_{i=1}^{k+1} p_i + \sum_{j=1}^{\ell+1} q_j \leq n. 
\eeq
\eenum
\end{Definition}
(We could of course equivalently say that \eref{nonover} is in the form \eref{dref3} as given in part (1) of Theorem \ref{mainfarkas-sdp}.)

Note that condition \eref{piqj} means that the identity blocks in the 
$(a_i)$ and $(y_j)$ sequences do not overlap.
The weakly infeasible instance in Example \ref{ex1} (when we choose  $\alpha=0$) 
is such an instance with $p_1 = q_1 = 1$ and $p_2 = q_2 = 0.$ 

A larger example follows: 
\bex \label{ex-nonoverlap} {\rm The SDP
\beq \label{ex-nonoverlap-sys}
\ba{rcl} 
a_i \bullet y & = & 0 \, (i=1,2), \\
a_3 \bullet y & = & -1, \\
            y & \succeq & 0    
\ena
\eeq
is weakly infeasible, where $a_1, a_2, a_3$ are given below:
\beq \label{ai-PQ}
a_1 \, = \, \bpx 1 & 0 & 0 & 0 & 0 \\
           0 & 0 & 0 & 0 & 0 \\
           0 & 0 & 0 & 0 & 0 \\
           0 & 0 & 0 & 0 & 0 \\
           0 & 0 & 0 & 0 & 0 \epx, \, a_2 \, = \, \bpx 5 & 1 & 2 & \un{2} & \un{0} \\
                                                       1 & 1 & 0 & 0 & 0 \\
						       2 & 0 & 0 & 0 & 0 \\
                                                       \un{2} & 0 & 0 & 0 & 0 \\
                                                       \un{0} & 0 & 0 & 0 & 0 \epx, \, 
a_3 = \bpx 3 & 2 & 1 & 3 &  - 2 \\ 2 & 0 & 0 & \un{0.5} & \un{1} \\ 1 & 0 & 1 & 0 & 0 \\ 3 & \un{0.5} & 0 & 0 & 0 \\ -2 & \un{1} & 0 & 0 & 0 \epx.
\eeq
(Some matrix entries are underlined, since we will return to this instance to explain our algorithm.)

Matrices $y_1, y_2, y_3$ (again with some underlined entries) 
below show that \eref{ex-nonoverlap-sys} is nonoverlapping weakly infeasible:
\beq \label{yj-PQ}
y_1 \, = \, \bpx 0 & 0 & 0 & 0 & 0 \\
           0 & 0 & 0 & 0 & 0 \\
           0 & 0 & 0 & 0 & 0 \\
           0 & 0 & 0 & 0 & 0 \\
           0 & 0 & 0 & 0 & 1 \epx, \, 
y_2 \, = \, \bpx 0 & 0 & 0 & 0 & \un{1} \\
           0 & 0 & 0 & 0 & \un{2} \\
           0 & 0 & 0 & 0 & 1 \\
           0 & 0 & 0 & 1 & 0 \\
           \un{1} & \un{2} & 1 & 0 & 0 \epx, \, 
y_3 \, = \, \bpx 0 & 0 & 0 & \un{0} & 3 \\
                 0 & 0 & 0 & \un{1} & 5 \\
                 0 & 0 & 0 & 4 & 1 \\
                 \un{0} & \un{1} & 4 & 1 & 2 \\
                 3 & 5 & 1 & 2 & 3 \epx.
\eeq
Note that here 
$$
k=\ell = 2, \, p_1 = p_2 = p_3 = 1, q_1 = q_2 = 1, \, q_3 = 0.
$$ 
It is of course easy to see directly that \eref{ex-nonoverlap-sys} is weakly infeasible.
However, if we generate matrices 
$a_i \, (i=4,5, \dots)$ such that 
$a_i \bullet y_1 = a_i \bullet y_2 = 0 \,$
and add the contraints $a_i \bullet y = a_i \bullet y_3 \,$ 
to \eref{ex-nonoverlap-sys}, then the resulting system is still weakly infeasible, since the first three constraints 
prove that it is infeasible, and $(y_1, y_2, y_3)$ prove that it is not strongly so.
However, weak infeasibility of the enlarged system is  difficult 
to confirm directly. (A simple dimension count shows that this way we can extend 
\eref{ex-nonoverlap-sys} to have $13$ constraints.)
}
\eex

To proceed with stating 
our algorithm, for $(a_1, \dots, a_{k+1}) \in \regfr(\psd{n}) \,$ with block sizes $p_1, \dots, p_{k+1}$ 
we denote the $i$th block containing $p_i$ integers by $P_i,$ i.e., 
$$
P_1 = \bigl\{1, \dots, p_1 \}, \, P_2 = \{ p_1 + 1, \dots, p_1 + p_2 \}, \dots
$$
For $(y_1, \dots, y_{\ell+1}) \in \revregfr(\psd{n})$ with block sizes $q_1, \dots, q_{\ell+1}$ 
we similarly denote the $j$th block containing $q_j$ integers by
$Q_j, \,$ i.e., 
$$
Q_1 =  \bigl\{ n - q_1 + 1, \dots, n \bigr\}, \, Q_2 = \bigl\{ n - q_1 - q_2 +1, \dots, n - q_1 \bigr\}, \dots
$$

For instance, in Example \ref{ex-nonoverlap}
\beq \label{piqinon}
P_1 = \{1\}, \, P_2 = \{2\}, \, P_3 = \{3\}, \, Q_1 = \{5\}, \, Q_2 = \{4\}, \, Q_3 = \emp.
\eeq
For $a \in \sym{n}$ and $P, Q \subseteq \{1, \dots, n \}$ we denote by 
$a(P,Q)$ the union of two blocks of $a:$ the first block is 
indexed by rows corresponding to $P$ and columns 
corresponding to $Q; $ and the second is the block symmetric with the first.
\co{ (i.e., it corresponds to rows indexed by $Q$ and columns indexed by $P$). 
For example, in $a_2$ in equation \eref{ai-PQ} the underlined entries correspond 
to $a_2(P_1, Q_1 \cup Q_2).$ }
We also write $a(P,:)$ to abbreviate $a(P, \{1, \dots, n \}).$

We are now ready to state our algorithm. 
Algorithm \ref{algo:non_overlapping} first chooses $(a_1, \dots, a_{k+1}) \in \regfr(\psd{n})$ 
and $(y_1, \dots, y_{\ell+1}) \in \revregfr(\psd{n})$ with block sizes $p_i$ and $q_j$ which satisfy inequality 
\eref{piqj}. 
Note that by Definitions   \ref{reg-def}  and \ref{revreg-def} we have 
$a_i(P_i,P_i) = I$ for all $i, \,$ and $y_j(Q_j,Q_j)=I$ for all $j, \,$ 
and apart from these blocks, only the entries in 
\beq \label{nonzero} 
a_i(P_1 \cup \dots \cup P_{i-1}, :) \, {\rm  and \,} y_j(Q_1 \cup \dots \cup Q_{j-1}, :)
\eeq
are allowed to be nonzero for all $i$ and $j.$

Algorithm \ref{algo:non_overlapping} fixes 
even the entries of all $a_i$ and of all $y_j \,$ in 
\eref{nonzero} to some arbitrary 
values, except for  the entries in 
\beq \label{freeai} 
a_i(P_{i-1}, Q_1 \cup \dots \cup Q_\ell) \,\, {\rm and} \,\, y_j(Q_{j-1}, P_1 \cup \dots \cup P_k) \, {\rm for \, all} \, i \geq 2, \, {\rm and \,} j \geq 2,
\eeq
which it leaves free. The algorithm then sets the entries in these free blocks to satisfy the equations   
\beq \label{aiyj}
\ba{rcl}
a_i \bullet y_j \, = \, \left\{ \ba{rl} 0 & {\rm if} \, (i,j) \neq (k+1, \ell+1), \\
                                      -1 & {\rm if} \, (i,j) = (k+1, \ell+1). 
                                \ena \right.
\ena
\eeq
This way we find the first $k+1$ equations in \eref{d} and the last part of the
algorithm generates the  remaining $m-k-1.$ 
\begin{algorithm}[H]
\caption{Nonoverlapping weakly infeasible SDP}
\label{algo:non_overlapping}
\begin{algorithmic}
\State  Choose integers $m, n, k, \ell, p_1, \dots, p_k, q_1, \dots, q_\ell > 0$ and $p_{k+1}, q_{\ell+1} \geq 0$ to satisfy \eref{piqj}. 
\State Let $(a_1, \dots, a_{k+1}) \in \regfr_{k+1}(\psd{n})$ with block sizes $p_1, \dots, p_{k+1}.$ 
\State Let $(y_1, \dots, y_{\ell+1}) \in \revregfr(\psd{n})$ with  block sizes $q_1, \dots, q_{\ell+1}.$ 
\State Fix the entries of the $a_i$ and of $y_j \,$ as described above, and leave the entries in \eref{freeai} free.
       \For {$i=2:(k+1)$} 
         \For {$j=2:(\ell+1)$}
               \State (*) Set $a_i(P_{i-1}, Q_{j-1})$ and $y_j(P_{i-1}, Q_{j-1})$
                to satisfy the equation  (\ref{aiyj}) for $a_i \bullet y_j.$
         \EndFor
      \EndFor
\State Find $a_{k+2}, \dots, a_m$ orthogonal to $y_1, \dots, y_\ell.$ 
\State Set $c = (0, \dots, 0, -1, a_{k+2} \bullet y_{\ell+1}, \dots, a_{m} \bullet y_{\ell+1})^T.$ 
\end{algorithmic}
\end{algorithm}

Algorithm \ref{algo:non_overlapping} can generate the weakly infeasible instance in 
Example \ref{ex1}  
by starting with only the offdiagonal entries of $a_2$ and $y_2$ free, then setting these 
to satisfy the equation $a_2 \bullet y_2 = -1.$ 

Algorithm \ref{algo:non_overlapping} can also generate the instance given in 
Example \ref{ex-nonoverlap}.
It starts with the 
underlined entries of the $a_i$ in (\ref{ai-PQ}) and of the $y_j$ in (\ref{yj-PQ}) 
free, and successively sets the entries in the following submatrices 
(note the definition of  $P_i$ and $Q_j$ in \eref{piqinon}) 
\benum
\item $a_2(P_1, Q_1)$ and $y_2(P_1, Q_1)$
\item $a_2(P_1, Q_2)$ and $y_3(P_1, Q_2)$
\item $a_3(P_2, Q_1)$ and $y_2(P_2, Q_1)$
\item $a_3(P_2, Q_2)$ and $y_3(P_2, Q_2)$
\eenum

\bth 
Algorithm \ref{algo:non_overlapping} always succeeds, and every nonoverlapping weakly infeasible instance 
is among its possible outputs.
\enth
\pf{} To show  that the algorithm always succeeds,  assume that at some point we execute Step (*). Since $P_{i-1}$ and $Q_{j-1}$ are nonempty, this step 
successfully satisfies equation (\ref{aiyj}). 
Let us consider an equation which involves $a_i, \,$ and which has been previously satisfied by the algorithm.
All such equations have left hand side 
$$
a_i \bullet y_t \, {\rm  with \,} t \leq j-1.
$$
Let us fix such a $t.$ We claim that 
\beq \label{ytzero} 
y_t(P_{i-1}, Q_{j-1}) = 0  
\eeq
holds. Indeed, first suppose $t=j-1.$ 
Then the only nonzero entries of $y_t$ corresponding  to columns indexed by $Q_{j-1}$ appear
in rows corresponding to $Q_1 \cup  \dots \cup Q_t$ (this is easiest to see by looking at the displayed formula in Definition \ref{revreg-def}). 
But inequality \eref{piqj} translates into $P_i \cap Q_j = \emp$ for all $i$ and $j$, so 
equation \eref{ytzero} follows. 

Next, assume $t < j-1.$ Then the only nonzero entries of $y_t$ corresponding  to columns indexed by $Q_{j-1}$ appear 
in rows corresponding to $Q_1 \cup  \dots \cup Q_{t-1}, \,$ so \eref{ytzero} follows similarly.

Thus all previously satisfied equations that contain $a_i$ remain satisfied.

It is trivial to prove that 
all nonoverlapping weakly infeasible instances are among the outputs: 
suppose that such an instance is identified by $(\bar{a}_1, \dots, \bar{a}_{m})$ and 
$(\bar{y}_1, \dots, \bar{y}_{\ell+1})$ with 
$(\bar{a}_1, \dots, \bar{a}_{k+1})$ having block sizes $p_1, \dots, p_{k+1}$, and 
the $\bar{y}_j$ having block sizes $q_j.$ 
Suppose that  Algorithm \ref{algo:non_overlapping} at the start sets all entries in $(a_1, \dots, a_{k+1})$ and 
$(y_1, \dots, y_{\ell+1})$ other than the ones in \eref{freeai} to the corresponding 
values in  the $(\bar{a}_i)$ and $(\bar{y}_j).$ Then there is a possible run of the algorithm 
which completes the $a_i$ and $y_j$ to be equal to the $\bar{a}_i$ and $\bar{y}_j.$ 
\qed

\section{Computational experiments} 
\label{sec-comp} 

To generate a test suite of challenging 
infeasible and weakly infeasible SDPs (in the dual form \eref{d}) we 
implemented  Algorithms \ref{algo:non_infeas_SDP} and \ref{algo:non_overlapping} in Matlab. 
We ran Algorithm \ref{algo:non_infeas_SDP}  with parameters 
\beq
n = 10, \, k=2, \, p_1 = 2, p_2 = 3, p_3 = 2, \, \, m = 10 \, {\rm or} \, m = 20,
\eeq
and  we call its outputs {\em infeasible} instances (these may be strongly or weakly infeasible). 
All entries in the generated instances are integers.

We ran Algorithm \ref{algo:non_overlapping} with parameters 
\beq
\begin{array}{ll}
n = 10, \, k=2, \, \ell = 1, \, p_1 = 2, \, p_2 = 3, \, p_3 =2, \, q_1 = 2, \, q_2 = 1 \\
m = 10 \; {\rm or} \; m = 20
\end{array}
\eeq
\co{
	\beq
n = 10, \, k=2, \, \ell = 1, \, p_1 = 2, \, p_2 = 3, \, p_3 =2, \, q_1 = 2, \, q_2 = 1, \,  m = 10 \, {\rm or} \, m = 20,
\eeq} 
and we call the instances it generates {\em weakly infeasible}. 
(These are guaranteed to be weakly infeasible.) 
We chose the components of $a_1, \dots, a_{k+1}$ in the  support of the
$y_j$ as integers in $[-2,2]$ so the entries in all $a_i$ 
turn out to be all integers.
The entries of the $y_j$ turn out to be 
rationals with denominators in the set $\{1, \dots, 5 \}.$ 

Hence one can verify the status of our instances in exact arithmetic.

To generate instances, in which the structure that leads to (weak) 
infeasibility is less readily apparent, we add the optional
\benum
\item[{\bf Messing step:}] Choose $t = (t_{ij}) \in \zad{m \times m}$ and $v = (v_{ij}) \in \zad{n \times n}$ 
random invertible matrices with entries in $[-2,2]$ and let 
$$
a_i \, = \, v^T \bigl( \sum_{j=1}^m t_{ij} a_j\bigr) v \; {\rm for} \;  i=1, \dots, m.
$$
\eenum
The $t$ matrix encodes elementary row operations performed on \eref{d}, and $v$ encodes a rotation. 

We call the instances output by Algorithms 2 and 3 {\em clean}, and 
the instances we find after the Messing step {\em messy.} 

The choices: ``clean/messy, infeasible/weakly infeasible,  $m=10/m=20$'' provide eight categories and 
we generated $100$ instances in each. 
We set the objective function as $I$ to ensure that the primal problem $(P)$ is feasible.

We tested four solvers: we first ran the solvers Sedumi, SDPT3 and MOSEK from the YALMIP environment, and 
the preprocessing algorithm of Permenter and Parrilo \cite{PerPar:14} interfaced with Sedumi. The latter is marked 
by ``PP+SEDUMI'' in our tables. 

As  the solvers consider our dual problem to be the primal, 
the only correct solution  status is ''primal infeasible.''
In  Tables \ref{table1} and \ref{table2} we report the number of solved instances out of $100$ for the various solvers. 

\renewcommand*{\arraystretch}{1.3}
\begin{table}[ht]
\begin{center}
\begin{tabular}{ |c|c|c|c|c|}
  \hline
  & \multicolumn{2}{|c|}{Infeasible} & \multicolumn{2}{c|}{Weakly Infeasible} \\
  \cline{2-5}
  & Clean & Messy & Clean & Messy \\
  \hline
  SEDUMI & 87 & 27& 0 & 0 \\ \hline
  SDPT3 & 10 & 5 & 0 & 0  \\ \hline
  MOSEK & 63 & 17 & 0 & 0\\ \hline
  PP+SEDUMI & 100 & 27 & 100 & 0\\
  \hline
\end{tabular}
\end{center}
\caption{Results with $n=10, m=10$}
\label{table1} 
\end{table}
\begin{table}[ht]
\begin{center}
\begin{tabular}{ |c|c|c|c|c|}
  \hline
  & \multicolumn{2}{c|}{Infeasible} & \multicolumn{2}{c|}{Weakly Infeasible} \\
  \cline{2-5}
  & Clean & Messy & Clean & Messy \\
  \hline
  SEDUMI & 100 & 100 & 1 & 0 \\ \hline
  SDPT3 & 100 & 96 & 0 & 0  \\ \hline
  MOSEK & 100 & 100 & 11 & 0\\ \hline
  PP+SEDUMI & 100 & 100 & 100 & 0\\
  \hline
\end{tabular}
\end{center}
\caption{Results with $n=10, m=20$} 
\label{table2} 
\end{table}

When the solvers do not report ''primal infeasible'', they mostly return an error status: 
for example, for the weakly infeasible, $n=10, \, m =20$ instances (both clean and messy) SDPT3
always reports an error status. 

We can see that 
\benum
\item The standalone solvers do better when $m$ goes from $10$ to $20$ as
for larger $m$ the portion of strongly infeasible instances is likely to be higher. 
\item The standalone solvers mostly fail on the weakly infeasible problems, 
though MOSEK detects infeasibility 
of some. These are ``almost'' strongly infeasible, i.e., the alternative system 
$(\dalt)$ is almost feasible. 
(Of course, in exact arithmetic $(\dalt)$  is infeasible.)

\item The preprocessing of \cite{PerPar:14} considerably helps Sedumi on the clean instances; it 
does not help, however, on the messy instances. 
\eenum
Clearly, a preprocessing algorithm like \cite{PerPar:14} could easily scan for entire facial reduction 
sequences in the input, and it is likely that some of the instances coming from applications also contain 
such sequences.
\co{
\section{Characterizing, and generating all weakly infeasible SDPs} 
\label{sec-generate}

\bth
The system \eref{d} is weakly infeasible, if and only if conditions (1) and (2) below hold:
\benum
\item it  has a reformulation 
\beq\label{dref2}\tag{$D_{\rm ref}$}
\ba{rcl}
\la a_{i}^{\prime}, y\ra & = & 0\,(i=1,\ldots,k)\\
\la a_{k+1}^{\prime}, y\ra & = & -1\\
\la a_{i}^{\prime}, y\ra & = & c_{i}^{\prime}\,(i=k+2,\ldots,m)\\
y & \geq_{K^*} & 0
\ena
\eeq
where $k \geq 0, \, (a_{1}^{\prime},\ldots,a_{k+1}^{\prime})\in\regfr(\psd{n}).$
\item there is $(y_1, \dots, y_{\ell+1}) \in \relregfr(\psd{n})$ such that 
\beq \label{condapyj} 
\ba{ccl}
A^{\prime *}y_{i} &=& 0 \,(i=1,\ldots,\ell)\\
A^{\prime *}y_{\ell+1} &=& c^\prime
\ena
\eeq
\eenum
Also, both facial reduction sequences can be chosen as pre-strict.
\enth

\pf{} The ``if'' direction is trivial.

To prove the ``Only if'' ditrection, we use 
Theorem \ref{mainfarkas-sdp} to choose $(A^\prime, c^\prime)$ to satisfy part (1). 
We then use 
Theorem \ref{mainfarkas} to choose 
$(y_1, \dots, y_{\ell+1})$ to satisfy the affine constraints of part (2), but with 
$(y_1, \dots, y_{\ell+1}) \in \fr(\psd{n}).$ 

In a general step we will have a reformulation with data $(A^\prime, c^\prime)$ and 
$(y_1, \dots, y_{\ell+1})$ which  satisfies the conditions of part (1) and the affine conditions of
part (2), with 
\beq  \label{condyjmod}
\ba{rcl}
(y_1, \dots, y_{\ell+1}) & \in & \fr_{\ell+1}(\psd{n}), \\ 
(y_1, \dots, y_{j}) & \in & \relregfr_{j}(\psd{n})
\ena
\eeq
for some $j \geq 0.$ 

We will call these conditions the {\em invariant conditions.} 

In $\ell+1$ steps we will transform this data 
so they will satisfy the invariant conditions with $j = \ell+1.$ 
To do so, we will repeatedly use the basic step 
\beq \label{transf} 
\ba{rcl} 
y_s & := & v^T y_s v  \; (s=1, \dots, \ell+1) \\
a_i & := & v^{-1} a_i^\prime v^{-T} \; (i=1, \dots, k+1).
\ena
\eeq
with various choices of a $v$ invertible matrix. 

We assume that $(a_1^\prime, \dots, a_{k+1}^\prime)$ are associated with sets $P_1, \dots, P_{k+1}, \,$ and 
$y_1, \dots, y_j$ with sets $Q_1, \dots, Q_{j}.$ 

\nin{\bf Step 1} We claim 
\beq
y_1(N \setminus P_{1:(k+1)}) \succeq 0, \, {\rm and \; the \; rest \; of \; } y_1 \; {\rm is \; zero}.
\eeq
Indeed, this follows from $y_1 \succeq 0, \,$ and $y_1$ being orthogonal to 
$a_1^\prime, \dots, a_{k+1}^\prime.$  

We let $q$ be an orthonormal basis of eigenvectors of 
$y_1(N \setminus P_{1:(k+1)})$ so the eigenvectors corresponding to zero eigenvalues come first, and 
$v = I \oplus q.$ We perform the basic step \eref{transf}, and 
we define $Q_1$ as the index set corresponding to positive diagonal components of $y_1.$ 
Clearly, then
\beq \label{y1q1} 
y_1(Q_1) \succ 0, \, {\rm diagonal, \, and \, the \, rest \, of \, } y_1 \, {\rm is \, zero},
\eeq
hence now \eref{condyjmod} holds with $j=1.$ 
For convenience, we define 
$$
P_{k+2} \, = \, N \setminus (P_{1:(k+1)}\cup Q_1).
$$
We refer to Figure \ref{y1before} to how $y_1$ looks before and after Step 1, assuming $k=1.$ 
\compilepdf{\begin{figure}
\centering
\includegraphics[scale=0.5]{ref1_revised}
\caption{The matrix $y_1$ before and after Step 1.}
\label{y1before}
\end{figure}
}

\nin{\bf Step 2} We claim 
\beqa \label{1st}
y_2(N \setminus Q_1) & \succeq & 0, \\ \label{2nd}
y_2(P_1) & = & 0, \\ \label{3rd}
y_2(P_1, N \setminus Q_1) & = & 0.
\eeqa
Indeed, \eref{1st} follows from \eref{y1q1}, and $y_2 \in (\psd{n} \cap y_1^\perp)^*.$ 
Equation \eref{2nd} follows from \eref{1st}, $a_1 \bullet y_2 = 0, \,$ and $P_1 \subseteq N \setminus Q_1.$ 
Finally, \eref{3rd} follows from \eref{1st} and \eref{2nd}.

Now for $t=2, \dots, k+1$ we perform the following operations: 
\benum
\item Define $q$ as the matrix of orthonormal eigenvectors of $y_2(P_t)$ 
with eigenvectors corresponding to positive eigenvalues coming first, and 
let $v = I \oplus q \oplus I.$ 
\item Apply the transformation \eref{transf}. 
\item \label{defv} Define $v$ as a full rank 
matrix so that 
\bit
\item[] left multiplying $y_2$ by $v^T$ adds multiples of the rows of  
$y_2(P_t, N)$ to the rows of 
$y_2(P_{t'}, N)$ to zero out the blocks $y_2(P_{t'}, P_t)$ for all $t < t' \leq k+2.$ 
\eit
\item Apply the transformation \eref{transf}. 
\eenum 
Let $Q_2$ be the index set corresponding to positive diagonal elements in 
$y_2(P_{2:(k+1)}).$   
Since in step \eref{defv}, 
left multiplying $a_i^\prime$ by $v^{-1}$ adds multiples  of the rows of 
$a_i^\prime$ corresponding to $P_{t^\prime}$ to the rows corresponding to 
$P_t$, the new $(a_1^\prime, \dots, a_{k+1}^\prime)$ is still in $\regfr(\psd{n}).$ 

So now the invariant conditions hold with $j=2.$ 

We refer to Figure \ref{y2before} to how $y_2$ looks before and after Step 2, assuming $k=1.$ 
\compilepdf{\begin{figure}
\centering
\includegraphics[scale=0.5]{ref2}
\caption{The matrix $y_2$ before and after Step 2: the index set $Q_2$ corresponds to the positive diagonal elements in 
$P_2$ and $P_3$} 
\label{y2before}
\end{figure}
}

\nin{\bf Step $j+1$, where $j \geq 2.$} We assume that 
the invariant conditions hold with $j. \,$ 
We claim 
\beqa \label{1stj}
y_{j+1}(N \setminus Q_{1:j}) & \succeq & 0, \\ \label{2ndj}
y_{j+1}(P_1) & = & 0, \\ \label{3rdj}
y_{j+1}(P_1, N \setminus Q_{1:j}) & = & 0.
\eeqa
Indeed, \eref{1stj} follows from 
$(y_1, \dots, y_{j}) \in \relregfr_{j}(\psd{n}),$ and 
$y_{j+1} \in (\psd{n} \cap y_1^\perp \cap \dots \cap y_{j}^\perp)^*.$ 
Equation \eref{2ndj} follows from \eref{1stj}, $a_1 \bullet y_{j+1} = 0, \,$ and 
$P_1 \subseteq N \setminus Q_{1:j}.$ 
Finally, \eref{3rdj} follows from \eref{1stj} and \eref{2ndj}.

Now, for $t=2, \dots, k+1$ we 
\bit
\item Define $q$ as the matrix of orthonormal eigenvectors of $y_{j+1}(P_t \setminus Q_{1:j})$ 
with eigenvectors corresponding to positive eigenvalues coming first, and 
let $v = I \oplus q \oplus I.$ 
\item Apply the transformation \eref{transf}. 
\item Define $v$ as a full rank 
matrix so that 
\bit
\item[] left multiplying $y_{j+1}$ by $v^T$ adds multiples of the rows of  
$y_{j+1}(P_t \setminus Q_{1:j}, N)$ to the rows of 
$y_{j+1}(P_{t^\prime} \setminus Q_{1:j}, N)$ to zero out the blocks 
$y_{j+1}(P_{t^\prime} \setminus Q_{1:j}, P_{t} \setminus Q_{1:j})$ for all $t < t' \leq k+2.$ 
\eit
\item Apply the transformation \eref{transf}. 
\eit 
Let $Q_{j+1}$ be the index set corresponding to positive diagonal elements in 
$y_{j+1}(P_{2:(k+2)}).$   
Clearly, now the invariant conditons hold with $j+1$ in place of $j.$

}

The SDP instances are available from 
\begin{verbatim}
www.unc.edu/~pataki/SDP.zip
\end{verbatim}

\section{Discussion and conclusion}
\label{sect-conclude}

Here we first mention how some of our results can be further extended.

First, the bound $\ell_K -1 \,$ in Theorem \ref{frlemma}, and 
the bounds in Corollary \ref{dimension}  have been 
improved by Lourenco et al in \cite{Lou:15}: they show that $\ell_K$ can be replaced 
by the length of the longest chain of faces which starts with $K$ and ends with a polyhedral face
of $K.$ 

Second,  we note that combining  Theorem \ref{mainfarkas} and part \eref{frkprop3} in Lemma \ref{frkprop}
we can write exact duals, and exact certificates of infeasibility
for more involved conic linear systems. For instance, we can prove that 
the system 
\beq
\ba{rcl} 
A_1 x & \leq_{K_1} & b_1, \\
A_2 x & \leq_{K_2} & b_2
\ena
\eeq
(where $K_1$ and $K_2$ are closed convex cones) is infeasible iff there is $k \geq 0$ and 
$(y_1, \dots, y_{k+1}) \in \fr_{k+1}(K_1)$ and $(z_1, \dots, z_{k+1}) \in \fr_{k+1}(K_2)$ such that 
$$
\ba{rccrcll}
A_1^*y_{i} + A_2^*z_{i} &=&0,&b_1^*y_{i} + b_2^*z_i &=&0\;(i=1,\ldots,k)\\
A_1^*y_{k+1} + A_2^*z_{k+1} &=&0,&b_1^*y_{k+1} + b_2^*z_{k+1} &=&-1.
\ena
$$
We next discuss how our machinery can be used to generate all infeasible 
conic LP instances over other classes of cones.
The general principle is: if the facial structure of $K$ is known
-- the case for all the cones 
over which we can efficiently optimize, e.g. for $\psd{n}$ -- then 
it is also easy to describe the corresponding facial reduction cone, and 
infeasible systems of the form \eref{dref} (in Theorems I and \ref{mainfarkas}). 
The facial structure is, in fact, trivial 
for a broad and useful class of cones, which we call {\em smooth} cones: we say that $K$ is a {\em smooth cone} if 
it is pointed, full-dimensional, and all faces distinct from $\{0\}$ and $K$ itself are one-dimensional
(i.e., extreme rays). Smooth cones include the $p$-order cones
$$
\{ \, (x_0, \bar{x}) \in \rad{} \ti \rad{n-1} \, | \, x_0 \geq \norm{\bar{x}}_p \, \},
$$
when $p \neq 1, + \infty, \,$ and the geometric and dual geometric cones
introduced by Glineur \cite{Glineur:01}. 

If $K$ is smooth, then any strict facial reduction sequence for $K$ is of length at most two. 
Thus, if $K$ is a direct product of 
smooth cones, then by  part (4) of Lemma \ref{frkprop} we can describe $\fr_k(K)$ 
and by part (1) of Theorem \ref{mainfarkas} 
we can easily generate all infeasible conic LP instances over $K.$ 

Consider the special case when $K$ is 
the direct product of second order cones, i.e.,
$$
K \, = \, K^* \, = \, \so{n_1}  \ti \dots \ti \so{n_q},
	$$
	where
	$$
	\so{n} \, = \,      \{ \, (x_0, \bar{x}) \in \rad{} \ti \rad{n-1} \, | \, x_0 \geq \norm{\bar{x}}_2 \, \}.   
	$$
	The cone $\so{n}$ has an interesting automorphism group. If 
	$x$ is on the boundary of $\so{n},$ then there is $T \in \aut(\so{n})$ such that
	$Tx = (1,1,0, \dots, 0)^T.$ 
	 If $x$ is in the interior of $\so{n},$ then there is $T \in \aut(\so{n})$ such that
	 $Tx = (1,0, \dots, 0)^T.$ We can define facial reduction sequences 
	 for the cone $K, \,$ which are in a sense ''regularized'' (just like we defined them 
	 for $\psd{n}$ in Definition \ref{reg-def}), and prove an analog of Theorem \ref{mainfarkas-sdp}: we plan to explore  this topic in a followup paper.
	 	
We can recover variants of Ramana's exact dual for SDP as follows. 
If $K$ is nice, then $(y_1, \dots, y_k) \in \fr_k(K)$ if and only if it has a decomposition 
\beq \label{frkrep} 
\begin{array}{rcl}
y_i & = & u_i + v_i \,  {\rm where} \\
(u_i,v_i) & \in & K^* \ti \tan(u_1 + \dots + u_{i-1}, K^*) \, {\rm for \, all \,} i. 
\end{array}
\eeq
Here $\tan(u, K^*)$ denotes the tangent space at $u \in K^*,$ and we set 
$v_0 = v_1 = 0.$  
Indeed, this decomposition follows from 
\cite[Section 4]{Pataki:13}. \cite{Pataki:13}. 
In turn, the cone $\psd{n}$ is nice, and 
\beq \label{tan-x}
\tan(u, \psd{n})  = \left\{ \, w + w^T \Bigm |  \, \bpx u    & w \\
w^T  & \beta I
\end{pmatrix} \,  \succeq 0   \,\, \text{for some} \;   
\beta  \in  \rad{} \right\},
\eeq
hence $\fr_k(\psd{n})$ is {\em semidefinite representable.} 


Thus the dual \eref{Dstrong} in Theorem \ref{sdlemma} 
provides a template to create exact duals of \eref{p} over known simple cones:
to write down such a dual, we only need to represent 
$\fr_k(K).$ Therefore it is natural to ask:
for what choices of $K$ is  
$\fr_k(K)$ representable in terms of $K^*, \, \psd{n}, \,$ or some other simple cones?

\begin{appendices}

\section{A definition of certificates}
\label{app-def} 

\co{An informal definition of {\em certificates} of 
certain properties of conic linear programs
was given in the Introduction,  and the results of the paper can be understood relying only on 
this informal definition. In this section we give a more 
rigorous definition of  certificates for completeness. A fully rigorous definition 
can be found in \cite{Cuckeretal:98}.}

In the Introduction we gave an informal definition of {\em certificates} of 
certain properties of conic linear programs, 
 and the results of the paper can be understood relying only on 
this informal definition. In this section we give a more 
rigorous definition of  certificates; 
a fully rigorous definition can be found in \cite{Cuckeretal:98}.

We define the set of primal instances as 
\beqa
\primal & = & \{ \, (A,b) \, | \, A: \rad{m} \rightarrow Y, \, b \in Y \, \}.
\eeqa
We assume that a map $A: \rad{m} \rightarrow Y$ is represented by a suitable matrix.

\begin{Definition}
Let 
$$
C: \primal \rightarrow Y'
$$
be a function with $Y'$ being a finite dimensional Euclidean space.
We say that $C$ 
provides an {\em exact certificate of infeasibility of \eref{p}}, if there is an algorithm,
a ''verifier'', which 
\begin{enumerate}
\item Takes as input $(A,b)$ and $C(A,b),$ where $(A,b) \in \primal;$
\item Outputs ''yes'' exactly when \eref{p} with data $(A,b)$ is infeasible; 
\item Takes a polynomial number of steps in the size of $(A,b)$ and $C(A,b).$
\end{enumerate}
\co{A "step" means either 
\begin{enumerate}
\item a usual arithmetic operation; or
\item checking membership in a 
set of the form $(K \cap L)^*, \,$ and 
$(K^* \cap L)^*, \,$ where $L$ and $J$ are subspaces. 
\end{enumerate}
}
A ''step'' means either a usual arithmetic operation; or checking membership in sets  of the form $(K \cap L)^*, \,$ and 
$(K^* \cap J)^*, \,$ where $L$ and $J$ are subspaces. 
By ''size'' of $a \in Y^\prime$ we mean 
the number of components of $a.$
\end{Definition}

We define an exact certificate of infeasibility of \eref{d}, and of other properties (say of weak infeasibility) of \eref{p} and of \eref{d} analogously.

The alert reader may notice that for our exact certificate of infeasibility given in  
part (3) of Theorem \ref{mainfarkas} it is enough to allow a 
verifier to check membership in a set of the form $(K \cap L)^*; $ 
and a similar restriction can be made for our exact certificate of 
infeasibility of \eref{d}.

\section{\hspace{-.2in}: Proof of Lemmas \ref{frkprop} and \ref{tclemma}}
\label{app-proofs}

\pf{of Lemma \ref{frkprop}} 

\pf{ of \eref{frkprop1}} It is clear  that $\fr_k(K)$ contains all nonnegative multiples of its elements,
so we only need to show that it is convex. To this end, we use the following Claim, 
whose proof is straightforward: 

\noindent{\bf Claim} If $C$ is a closed, convex cone and $y, z \in C^*,$ then 
$$
C \cap (y+z)^\perp = C \cap y^\perp \cap z^\perp.
$$
\qed 

Let $(y_1, \dots, y_k), (z_1, \dots, z_k) \in \fr_k(K).$ We will prove
\beq \label{elaine2} 
(y_1 + z_1, \dots, y_k + z_k) \in \fr_k(K),
\eeq
which will clearly imply \eref{frkprop1}. To start, for brevity, for $i=0, \dots, k$ 
we set 
\beqast
K_{y,i} & = & K \cap y_1^\perp \cap \dots \cap y_i^\perp, \\
K_{z,i} & = & K \cap z_1^\perp \cap \dots \cap z_i^\perp, \\
K_{y+z,i} & = & K \cap (y_1+z_1)^\perp \cap \dots \cap (y_i+z_i)^\perp
\eeqast
(with the understanding that all these cones equal $K$ when $i=0$). 
We first prove the relations 
\beqa \label{yiK}
y_i & \in &  K_{y+z,i-1}^*, \\ \label{ziK} 
z_i & \in &  K_{y+z,i-1}^*, \\ \label{Ky1z1}
K_{y+z,i} & = & K_{y,i} \cap K_{z,i}
\eeqa
for $i=1, \dots, k.$ For $i=1$ the first two hold by definition,
and \eref{Ky1z1} follows from the Claim. Suppose now that $i \geq 2$ and (\ref{yiK}) through (\ref{Ky1z1}) hold
with $i-1$ in place of $i.$
Then 
$$
y_i \, \in \, K_{y,i-1}^* \, \subseteq \, (K_{y,i-1} \cap K_{z,i-1})^* \, =  \, K_{y+z,i-1}^*,
$$
where the first containment is by definition, the inclusion is trivial, and the equality is by using the 
induction hypothesis. This proves \eref{yiK} and 
equation \eref{ziK} holds analogously. 

Hence
\beqast
K_{y+z,i} & = & K_{y+z,i-1} \cap (y_i + z_i)^\perp \\
          & = & K_{y+z,i-1}    \cap y_i^\perp \cap z_i^\perp \\ 
          & = & K_{y,i-1} \cap K_{z,i-1} \cap y_i^\perp \cap z_i^\perp \\
          & = & K_{y,i} \cap K_{z,i},
\eeqast
where the first equation is by definition. The second follows 
since by \eref{yiK} and \eref{ziK}, and since $K_{y+z,i-1}$ is a closed convex cone,
 we can use the Claim with $C = K_{y+z,i-1}, \, y = y_i, \, z = z_i. \,$
The third is by the induction hypothesis, and the last is by definition. 
This completes the proof of \eref{Ky1z1}. 

Now by \eref{yiK}, \eref{ziK}  and since $K_{y+z,i-1}^*$ is a convex cone, we deduce that 
$$
y_i + z_i \in K_{y+z,i-1}^* \,\, {\rm holds \, for \,} i=1, \dots, k.
$$
This proves \eref{elaine2}, and completes the proof of (1).

\pf{of \eref{frkprop15}} 
Let $L = K \cap -K,$ assume that $K$ is not a subspace, i.e., 
$K \neq L, $ and also assume $k \geq 2.$ 
Let us choose a sequence 
$\{ y_{1i} \} \subseteq \ri K^*, \,$ s.t. $y_{1i} \rightarrow 0.$ Then 
$$
\ba{rclrrcl}
K \cap y_{1i}^\perp & = & L & \Rightarrow & (K \cap y_{1i}^\perp)^* & = & L^\perp \,\,  {\rm for \, all} \; i,  \\
K \cap 0^\perp & = & K & \Rightarrow & (K \cap 0^\perp)^* & = & K^*. 
\ena
$$
Let $y_2 \in L^\perp \setminus K^*.$ (Such a $y_2$ exists, since $K^* \neq L^\perp.$) 
Then $(y_{i1}, y_2, 0, \dots, 0) \in \fr_k(K), \,$ and it converges to 
$(0, y_2, 0, \dots, 0) \not\in \fr_k(K). \,$

Conversely, if $K$ is a subspace, then an easy calculation shows 
that so is $\fr_k(K), \,$ which is hence closed. 

\pf{ of \eref{frkprop2}} Let us fix $T \in \aut(K)$ and let $S$ be an arbitary set. Then we claim that 
\beqa \label{star}  
T^{-1}S^* & = & (T^*S)^*  \\ \label{perp} 
T^{-1}S^\perp & = & (T^*S)^\perp \\ \label{last} 
T^* (K \cap S^\perp)^* & = & (K \cap (T^*S)^\perp)^*  
\eeqa
hold. Statement \eref{star} follows, since  
\beq \nonumber 
\ba{rcll}
 y \in T^{-1} S^* & \LRA & Ty \in S^*  \\
& \LRA & \la T y, x \ra \geq 0 & \forall x \in S \,  \\
& \LRA & \la y, T^*x \ra \geq 0 & \forall x \in S \\ 
& \LRA & y \in (T^*S)^*, 
\ena
\eeq
and \eref{perp}  follows analogously. Statement \eref{last} follows by 
\beq \nonumber 
\ba{rcll} 
T^* (K \cap S^\perp)^* & = & ( T^{-1}(K \cap S^\perp))^* \\
                                   & = & ( T^{-1}K \cap T^{-1}S^\perp)^* \\
                                   & = & (K \cap T^{-1}S^\perp)^*  \\
                                   & = & (K \cap (T^*S)^\perp)^*, 
\ena
\eeq
where 
in the first equation we used (\ref{star}) with $T^{-*}$ in place of $T$ and $K \cap S^\perp$ in place of $S$.  The seond equation is trivial, and in the third we used 
$T^{-1}K = K.$ In the last we used (\ref{perp}). 

Now let $(y_1, \dots, y_k) \in \fr_k(K),$ 
and  $S_i = \{y_1, \dots, y_{i-1}\}$ for $i = 1, \dots, k.$ Then by definition we have 
$
y_i \, \in \, (K \cap S_i^\perp)^* \, {\rm for \, all \,} i. 
$
Hence for all $i$ we have 
$$
T^* y_i \in T^*(K \cap S_i^\perp)^* = (K \cap (T^*S_i)^\perp )^*,
$$
where the equation follows from \eref{last}. Thus 
$
(T^*y_1, \dots, T^*y_k) \in \fr_k(K),
$
as required. 

\co{
\pf{ of \eref{frkprop3}} 
The equivalence of \eref{KC} and of \eref{KC2} holds for $k=1, \,$ so let us asssume $k \geq 2$ and that we proved it with $1, \dots, k-1$ in place of $k.$ 

Statement  \eref{KC} is equivalent to 
\beq\label{second-1}
\bigl( (y_1, z_1), \dots, (y_{k-1}, z_{k-1})\bigr) \, \in \, FR_{k-1}(K \times C)
\eeq
and
\beq \label{second}
(y_{k}, z_{k}) \, \in \, \bigl((K \ti C) \cap (y_1, z_1)^\perp \cap \dots \cap (y_{k-1}, z_{k-1})^\perp \bigr)^*.
\eeq
However, 
\beqast
(K \ti C) \cap (y_1, z_1)^\perp \cap \dots \cap (y_{k-1}, z_{k-1})^\perp & = & \bigl((K \ti C) \cap (y_1, z_1)^\perp \cap \dots \cap (y_{k-2}, z_{k-2})^\perp \bigr) \cap (y_{k-1}, z_{k-1})^\perp \\
        & = & 
\eeqast

By the inductive hypothesis the set on the right hand side of \eref{second} is 
$$
(K \cap y_1^\perp \cap \dots \cap y_{k}^\perp)^* \times (C \cap z_1^\perp \cap \dots \cap z_{k}^\perp)^*,
$$
and this completes the proof.
}

\pf{ of \eref{frkprop3}} 
In this proof, for brevity,  we will use the notation
\beqast
K_{y,i} & = & K \cap y_1^\perp \cap \dots \cap y_i^\perp, \\
C_{z,i} & = & C \cap z_1^\perp \cap \dots \cap z_i^\perp, \\
(K \ti C)_{y,z,i} & = & (K \ti C)  \cap (y_1, z_1)^\perp \cap \dots \cap (y_{i}, z_{i})^\perp 
\eeqast
for all $i \geq 0$ (with the understanding that these sets equal $K, C, \,$ and $K \times C, \,$ respectively, when $i=0$). 

We will prove the equivalence $\eref{KC} \LRA  \eref{KC2}$ together with the relation 
\beq \label{KC3} 
(K \ti C)_{y,z,i} \, = \, K_{y,i} \times C_{z,i} \, {\rm for} \, i \leq k-1.
\eeq

Clearly, both hold for $k=1, \,$ so let us assume $k \geq 2$ and that we proved them with $k-1$ 
in place of $k.$ By definition,  \eref{KC} is equivalent to 
\beq \label{KC-1}
\bigl( (y_1, z_1), \dots, (y_{k-1}, z_{k-1}) \bigr) \in \fr_{k-1}(K \ti C)
\eeq
and 
\beq \label{KC-2}
(y_k, z_k) \in (K \ti C)_{y,z,k-1}^*. 
\eeq
By the induction hypothesis \eref{KC-1} is equivalent to 
$(y_1, \dots, y_{k-1}) \in \fr_{k-1}(K)$ and $(z_1, \dots, z_{k-1}) \in \fr_{k-1}(C).$
		So the proof is complete, if we show 
	\beq \label{KC4} 
	(K \ti C)_{y,z,k-1} \, = \, K_{y,k-1} \times C_{z,k-1}. 
	\eeq
	To prove (\ref{KC4}) we see that  
		\beqast
		(K \ti C)_{y,z,k-1}  & = &  	(K \ti C)_{y,z,k-2} \cap (y_{k-1}, z_{k-1})^\perp \\
		                            & = &   (K_{y, k-2} \ti  C_{z,k-2})  \cap (y_{k-1}, z_{k-1})^\perp \\
		                            & = & (K_{y, k-2} \cap y_{k-1}^\perp) \ti (K_{z, k-2} \cap z_{k-1}^\perp) \\
		                            & = & K_{y,k-1} \times C_{z,k-1},
		                            \eeqast
	where the first equation is by definition, the second  comes from the induction hypothesis, and the third follows from 
	$y_{k-1} \in K_{y,k-2}^*, \, z_{k-1} \in C_{z,k-2}^*.$ Thus the proof is complete.

\pf{of Lemma \ref{tclemma}} First let us note that $T \in \aut(\psd{n})$ (cf. equation \eref{def-autk}) 
iff $T(x) = t^T x t$ for some invertible matrix $t.$

We prove the lemma by induction. Suppose that $\ell \geq 0$ is an integer, and 
we computed a $t$ invertible matrix such that 
\beqa \label{tytk}
(t^Ty_1t, \dots, t^Ty_{k}t) & \in & \fr(\psd{n}), \\ \label{tytl}
(t^Ty_1t, \dots, t^Ty_{\ell}t) & \in & \regfr(\psd{n}),
\eeqa
and the block sizes in the latter sequence are 
$p_1, \dots, p_{\ell}, \,$ respectively. Both of these statements hold with $\ell =0.$ 
 If $\ell = k, \,$ we stop.

Otherwise, define $p  := p_1 + \dots + p_{\ell}$ and $y_i^\prime := t^Ty_it$ for $i=1, \dots, k.$ 
Let 
$$K = \psd{n} \cap y_1^{\prime \perp} \cap \dots \cap y_{\ell}^{\prime \perp}.
$$ 
Then $y_{\ell+1}' \in K^*, \,$ and $K$ and $K^*$ are of the form 
$$
K = \bordermatrix{& p  &   n-p   \cr
                  & 0  & 0       \cr
                  & 0  & \oplus \cr}, \, K^* = \bordermatrix{& p  &   n-p   \cr
                                                             & \ti  & \ti       \cr
                                                             & \ti  & \oplus \cr},
$$
where, again, the symbol $\oplus$ stands for a psd submatrix, and $\ti$ for a submatrix with arbitrary elements.

Let $z$ be the lower $n-p$ by $n-p$ block of $y_{\ell+1}^\prime.$ Since $z \succeq 0,$ 
there is 
a $q$ invertible matrix such that 
$$
q^T z q \, = \, \bpx I_{p_{\ell+1}} & 0 \\
                     0              & 0 \epx,
$$
where $p_{\ell+1}$ is the rank of $z.$ 

Let $v := I_p \oplus q$ and replace $t$ by $tv.$ Then by part \eref{frkprop2} in 
Lemma  \ref{frkprop} statement \eref{tytk} still holds, and by the choice of $v$ equation \eref{tytl} now 
holds with $\ell+1$ in place  of $\ell.$ This completes the proof. 
\qed

\section{\hspace{-.2in}: Proof of Theorem \ref{nonclosed}}
\label{app-nonclosed}
\pf{of (1)}  Let us assume that condition \eref{radir} is violated;  
we  will construct $(a_1, a_2) \in \fr(K^*)$ and 
$(y_1, \dots, y_{\ell+1}) \in \fr(K)$ that satisfy 
(\ref{aiyj-1}) and (\ref{aiyj-2}) (with $k=1, \,$ and 
$\ell$ equal to the degree of singularity of 	$\R(A) \cap K$).  

First, we choose  
\beqast
a_1 & \in & \ri  (\R(A) \cap K), \\
a_2 & \in & \R(A) \cap (\cl \dir(a_1, K) \setminus \dir(a_1,K)),
\eeqast
and let $F$ be the minimal cone  of $\R(A) \cap K \,$ (i.e., the smallest face of $K$ that contains 
$a_1$). Then 
$$
(K^* \cap a_1^\perp)^* \, = \, (K^* \cap F^\perp)^* \, = \, \cl \dir(a_1, K),
$$
where the first equality comes from $a_1 \in \ri F$ and the second can be found e.g., in 
\cite{Pataki:00A}. Hence 
\beqast
a_1, a_2 & \in & \R(A),  \\
(a_1, a_2) & \in & \fr(K^*)
\eeqast
hold. We next choose the $y_j.$ First we pick 
$(y_1, \dots, y_\ell) \in \fr(K)$  such that $y_1, \dots, y_\ell \in \N(A^*)$ and 
$$
F = K \cap y_1^\perp \cap \dots \cap y_\ell^\perp.
$$
Since $a_2 \not \in \lin F \,$ 
(otherwise $a_2$ would be in $\dir(a_1,K)$)
we can then choose $y_{\ell+1} \in F^\perp$ such that 
\beqast
\la a_1, y_{\ell+1} \ra & = & 0, \\
\la a_2, y_{\ell+1} \ra & = & -1 
\eeqast
hold. Thus $(a_1, a_2)$ and $(y_1, \dots, y_{\ell+1})$ are as required, and
the proof is complete.

\pf{of (2)}  We fix $F$ and $G$ as stated.
Since $G$ is not exposed, and $F$ is the smallest exposed face of $K$ that contains $G,$ we have
$$
G \subsetneq F, \, K^* \cap G^\perp = K^* \cap F^\perp
$$
(see equation (\ref{FGperp}) and the discussion afterwards). 
For brevity, let us define $F^\tri = K^* \cap F^\perp, \, G^\tri = K^* \cap G^\perp,$ 
and for a face $H$ of $K^*$ we define 
$H^\tri = K \cap H^\perp.$ 
Thus, since $F$ is an exposed face, we also have 
$$
F^{\tri \tri} = F. 
$$
We will choose $(a_1, a_2)$ and $(y_1, y_2)$ such that
\beqa 
\label{frk*2}
(a_1, a_2)& \in & \fr(K^*), \\ 
\label{riG} 
a_1, a_2     & \in & \lin F, \\
\label{yfrk} 
(y_1, y_2) & \in & \fr(K), \\
\label{y1Fperp} 
y_1            & \in & F^\perp, \\ 
\label{a1y2} 
\la a_1, y_2 \ra & = & 0, \\
\label{a2y2}
\la a_2, y_2 \ra & = & -1
\eeqa
(i.e., to satisfy (\ref{aiyj-nice-1}) and (\ref{aiyj-nice-2}) with $k=1, \, \ell = 1$).  

We first choose $y_1 \in \ri F^\tri. \,$ Next, since $G \subsetneq F, \,$ we can choose  
$y_2 \, \in \, (F^* \cap G^\perp) \setminus F^\perp.$ Hence 
$$
K \cap y_1^\perp = K \cap (F^\tri)^\perp = F^{\tri \tri} = F,
$$
where the first equation comes from $y_1 \in \ri F^\tri. \,$  
We thus satisfied  (\ref{yfrk}) and (\ref{y1Fperp}). 

Next we choose $a_1$ and $a_2:$ 
we choose $a_1 \in \ri G, \,$ and $a_2 \in \lin F$ to satisfy \eref{a2y2} (this can be done since
$y_2 \not \in F^\perp$). Thus (\ref{riG}) and \eref{a1y2} also  hold. 
We  claim that \eref{frk*2} holds as well. 
To see this, we observe 
$$
K^* \cap a_1^\perp \, = \, K^* \cap G^\perp \, = \, K^* \cap F^\perp,
$$
where the first equation follows from $a_1 \in \ri G, \,$ and the second from $F^\tri = G^\tri.$ Hence 
$$
(K^* \cap a_1^\perp)^* \, = \, (K^* \cap F^\perp)^* \, \supseteq \, \lin F \, \ni a_2, 
$$
so \eref{frk*2} follows, and this completes the proof. 
\qed

\end{appendices}

{\bf Acknowledgements}  
We are grateful to the referees, the Associate Editor, and Melody Zhu  
for their insightful comments, and to 
Imre P\'{o}lik for his help in our work with the SDP solvers.

\bibliography{mysdp}

\begin{thebibliography}{10}

\bibitem{Aus:96}
Alfred Auslender.
\newblock Closedness criteria for the image of a closed set by a linear
  operator.
\newblock {\em Numer. Funct. Anal. Optim.}, 17:503--515, 1996.

\bibitem{BarCar:75}
George~Phillip Barker and David Carlson.
\newblock Cones of diagonally dominant matrices.
\newblock {\em Pacific J. Math.}, 57:15--32, 1975.

\bibitem{BausBor:99}
Heinz Bauschke and Jonathan~M. Borwein.
\newblock Conical open mapping theorems and regularity.
\newblock In {\em Proceedings of the Centre for Mathematics and its
  Applications 36}, pages 1--10. Australian National University, 1999.

\bibitem{Berman:73}
Abraham Berman.
\newblock {\em Cones, Matrices and Mathematical Programming}.
\newblock Springer-Verlag, Berlin, New York, 1973.

\bibitem{BertTseng:07}
Dimitri Bertsekas and Paul Tseng.
\newblock Set intersection theorems and existence of optimal solutions.
\newblock {\em Math. Program.}, 110:287--314, 2007.

\bibitem{Cuckeretal:98}
Lenore Blum, Felipe Cucker, Michael Shub, and Stephen Smale.
\newblock {\em Complexity and Real Computation}.
\newblock Springer, 1998.

\bibitem{BonnShap:00}
Fr\'ed\'eric~J. Bonnans and Alexander Shapiro.
\newblock {\em Perturbation analysis of optimization problems}.
\newblock Springer Series in Operations Research. Springer-Verlag, 2000.

\bibitem{BorLewis:00}
Jonathan~M. Borwein and Adrian~S. Lewis.
\newblock {\em Convex Analysis and Nonlinear Optimization: Theory and
  Examples}.
\newblock CMS Books in Mathematics. Springer, 2000.

\bibitem{BorweinMoors:09}
Jonathan~M. Borwein and Warren~B. Moors.
\newblock Stability of closedness of convex cones under linear mappings.
\newblock {\em J. Convex Anal.}, 16(3--4):699--705, 2009.

\bibitem{BorweinMoors:10}
Jonathan~M. Borwein and Warren~B. Moors.
\newblock Stability of closedness of convex cones under linear mappings ii.
\newblock {\em Journal of Nonlinear Analysis and Optimization: Theory \&
  Applications}, 1(1), 2010.

\bibitem{BorWolk:81B}
Jonathan~M. Borwein and Henry Wolkowicz.
\newblock Facial reduction for a cone-convex programming problem.
\newblock {\em J. Aust. Math. Soc.}, 30:369--380, 1981.

\bibitem{BorWolk:81}
Jonathan~M. Borwein and Henry Wolkowicz.
\newblock Regularizing the abstract convex program.
\newblock {\em J. Math. Anal. App.}, 83:495--530, 1981.

\bibitem{CheWolkSchurr:12}
Vris Cheung, Henry Wolkowicz, and Simon Schurr.
\newblock Preprocessing and regularization for degenerate semidefinite
  programs.
\newblock In David Bailey, Heinz~H. Bauschke, Frank Garvan, Michel Th\'era,
  Jon~D. Vanderwerff, and Henry Wolkowicz, editors, {\em Proceedings of
  Jonfest: a conference in honour of the 60th birthday of Jon Borwein}.
  Springer, 2013.

\bibitem{ChuaTuncel:08}
Check-Beng Chua and Levent Tun\c{c}el.
\newblock Invariance and efficiency of convex representations.
\newblock {\em Math. Program. B}, 111:113--140, 2008.

\bibitem{Dima:15}
Dimitry Drusviyatsky, G\'abor Pataki, and Henry Wolkowicz.
\newblock Coordinate shadows of semi-definite and euclidean distance matrices.
\newblock {\em SIAM J. Opt.}, 25(2):1160--1178, 2015.

\bibitem{Guler:10}
Osman G$\ddot{\mathrm{u}}$ler.
\newblock {\em Foundations of Optimization}.
\newblock Graduate Texts in Mathematics. Springer, 2010.

\bibitem{Glineur:01}
Francois Glineur.
\newblock Proving strong duality for geometric optimization using a conic
  formulation.
\newblock {\em Ann. Oper. Res.}, 105(2):155--184, 2001.

\bibitem{GortlerThurston:14}
Steven~J. Gortler and Dylan~P. Thurston.
\newblock Characterizing the universal rigidity of generic frameworks.
\newblock {\em Discrete Comput. Geometry}, 51(4), 2014.

\bibitem{KlepSchw:12}
Igor Klep and Markus Schweighofer.
\newblock An exact duality theory for semidefinite programming based on sums of
  squares.
\newblock {\em Math. Oper. Res.}, 38(3):569--590, 2013.

\bibitem{KrisWolk:10}
Nathan Krislock and Henry Wolkowicz.
\newblock Explicit sensor network localization using semidefinite
  representations and facial reductions.
\newblock {\em SIAM J. Opt.}, 20:2679--2708, 2010.

\bibitem{LiuPataki:15}
Minghui Liu and Gábor Pataki.
\newblock Exact duality in semidefinite programming based on elementary
  reformulations.
\newblock {\em SIAM J. Opt.}, 25(3):1441--1454, 2015.

\bibitem{Lou:15}
Bruno Lourenco, Masakazu Muramatsu, and Takashi Tsuchiya.
\newblock Facial reduction and partial polyhedrality.
\newblock {\em Optimization Online}, 2015.

\bibitem{Lourenco:13}
Bruno Lourenco, Masakazu Muramatsu, and Takashi Tsuchiya.
\newblock A structural geometrical analysis of weakly infeasible {SDPs}.
\newblock {\em Journal of the Operations Research Society of Japan},
  59(3):241--257, 2015.

\bibitem{Pataki:00A}
G\'abor Pataki.
\newblock The geometry of semidefinite programming.
\newblock In Romesh Saigal, Lieven Vandenberghe, and Henry Wolkowicz, editors,
  {\em Handbook of semidefinite programming}. Kluwer Academic Publishers, also
  available from www.unc.edu/\textasciitilde pataki, 2000.

\bibitem{Pataki:07}
G{\'a}bor Pataki.
\newblock On the closedness of the linear image of a closed convex cone.
\newblock {\em Math. Oper. Res.}, 32(2):395--412, 2007.

\bibitem{Pataki:12}
G\'abor Pataki.
\newblock On the connection of facially exposed and nice cones.
\newblock {\em J. Math. Anal. App.}, 400:211--221, 2013.

\bibitem{Pataki:13}
G\'abor Pataki.
\newblock Strong duality in conic linear programming: facial reduction and
  extended duals.
\newblock In David Bailey, Heinz~H. Bauschke, Frank Garvan, Michel Th\'era,
  Jon~D. Vanderwerff, and Henry Wolkowicz, editors, {\em Proceedings of
  Jonfest: a conference in honour of the 60th birthday of Jon Borwein}.
  Springer, also available from http://arxiv.org/abs/1301.7717, 2013.

\bibitem{Pataki:17}
Gábor Pataki.
\newblock Bad semidefinite programs: they all look the same.
\newblock {\em SIAM J. Opt.}, 27(1):146--172, 2017.

\bibitem{PerPar:14}
Frank Permenter and Pablo Parrilo.
\newblock Partial facial reduction: simplified, equivalent sdps via
  approximations of the psd cone.
\newblock Technical report, http://arxiv.org/abs/1408.4685, 2014.

\bibitem{PolikTerlaky:09}
Imre P\'olik and Tam\'as Terlaky.
\newblock Exact duality for optimization over symmetric cones.
\newblock Technical report, Lehigh University, Betlehem, PA, USA, 2009.

\bibitem{ProvanShier:96}
J.~Scott Provan and Douglas~R. Shier.
\newblock A paradigm for listing (s, t)-cuts in graphs.
\newblock {\em Algorithmica}, 15(4):351--372, 1996.

\bibitem{Ramana:97}
Motakuri~V. Ramana.
\newblock An exact duality theory for semidefinite programming and its
  complexity implications.
\newblock {\em Math. Program. Ser. B}, 77:129--162, 1997.

\bibitem{RaFreund:96}
Motakuri~V. Ramana and Robert Freund.
\newblock On the elsd duality theory for sdp.
\newblock Technical report, MIT, 1996.

\bibitem{RaTuWo:97}
Motakuri~V. Ramana, Levent Tun\c{c}el, and Henry Wolkowicz.
\newblock Strong duality for semidefinite programming.
\newblock {\em SIAM J. Opt.}, 7(3):641--662, 1997.

\bibitem{ReadTarj:75}
R.C. Read and R.E. Tarjan.
\newblock Bounds on backtrack algorithms for listing cycles, paths, and
  spanning trees.
\newblock {\em Networks}, 5:237--252, 1975.

\bibitem{Ren:01}
James Renegar.
\newblock {\em A Mathematical View of Interior-Point Methods in Convex
  Optimization}.
\newblock MPS-SIAM Series on Optimization. SIAM, Philadelphia, USA, 2001.

\bibitem{Rockafellar:70}
Tyrrel~R. Rockafellar.
\newblock {\em Convex Analysis}.
\newblock Princeton University Press, Princeton, NJ, USA, 1970.

\bibitem{Vera:13}
Vera Roshchina.
\newblock Facially exposed cones are not nice in general.
\newblock {\em SIAM J. Opt.}, 24:257--268, 2014.

\bibitem{Waki:12}
Hayato Waki.
\newblock How to generate weakly infeasible semidefinite programs via
  {L}asserre's relaxations for polynomial optimization.
\newblock {\em Optim. Lett.}, 6(8):1883--1896, 2012.

\bibitem{WakiMura:12}
Hayato Waki and Masakazu Muramatsu.
\newblock Facial reduction algorithms for conic optimization problems.
\newblock {\em J. Optim. Theory Appl.}, 158(1):188--215, 2013.

\end{thebibliography}

\end{document}